\documentclass[twocolumn]{autart}    
 \usepackage{mathrsfs}
\usepackage{amsfonts}
\usepackage{stmaryrd}
\usepackage{amsmath, amssymb}
\usepackage{cases}                                    
\usepackage{appendix}                                     

\usepackage{colortbl}
\usepackage{graphicx}                                     

\usepackage{graphicx}          

\begin{document}

\begin{frontmatter}

\title{Optimal Control for Discrete-time NCSs with Input Delay and Markovian Packet Losses: Hold-Input Case\thanksref{footnoteinfo}} 

\thanks[footnoteinfo]{This work is supported by Research Grants Council of Hong Kong under grant 15215319, the National Natural Science Foundation of China under Grants 61633014, U1701264 and the foundation for innovative research groups of national natural science foundation of China (61821004). Corresponding author Huanshui Zhang.}

\author[Paestum]{Hongdan Li}\ead{lhd200908@163.com},    
\author[Rome]{Xun Li}\ead{malixun@polyu.edu.hk.},               
\author[Paestum]{Huanshui Zhang}\ead{hszhang@sdu.edu.cn}  

\address[Paestum]{School of Control Science and Engineering,
Shandong University, Jinan, Shandong, P.R.China 250061.}
\address[Rome]{Department of Applied Mathematics, The Hong Kong Polytechnic University,
Kowloon, Hong Kong, P.R. China}         


\begin{abstract}                          
This paper is concerned with the linear quadratic optimal control problem for networked system simultaneously with input delay and Markovian dropout. Different from the results in the literature, we consider the hold-input strategy, which is much more computationally complicated than zero-input strategy, but much better in most cases especially in the transition phase. Necessary and sufficient conditions for the solvability of optimal control problem over a finite horizon are given by the coupled difference Riccati-type equations. Moreover, the networked control system is mean-square stability if and only if the coupled algebraic Riccati-type equations have a particular solution.  The key technique in this paper is to tackle the forward and backward difference equations, which are more difficult to be dealt with, due to the adaptability of controller and the temporal correlation caused by simultaneous input delay and Markovian jump.
\end{abstract}

\end{frontmatter}

\section{Introduction}
A network control system (NCS, for short), known as communication and control system, is a fully distributed and networked real-time feedback control system. Since the concept of NCSs was proposed in the early 1990s, it has attracted attention. For example, Guo et al \cite{14} discussed the networked control problems for linear discrete systems whose network mediums for the actuators are constrained. Using the Lyapunov-Krasovskii functional method, Yue et al \cite{3} considered the disturbance attenuation problem for NCSs. At the same time, it raised new challenges to traditional control system theory and applications. In an NCS, multiple network nodes share a network channel. Due to the limited network bandwidth and irregular changes in data traffic in the network, data collisions and network congestion often occur when multiple nodes exchange data through the network. Therefore, packet dropout and time delay will inevitably occur.\\
Actually, there are two different kinds of compensation strategies considered in the literature: the zero-input (i.e., zero value is directly adopted by the actuator input) and the hold-input (i.e., the latest available control signal stored in the actuator buffer is used). For zero-input case, Imer et al \cite{Imer: 06} discussed the optimal control problem for linear system with packet dropout under TCP and UDP protocols. Liang and Xu \cite{6} focused on the optimal control problems for NCSs which are simultaneously controlled by the remote controller and the local controller. Sufficient conditions for stability of network communication models with packet dropout were studied by  Montestruque and Antsaklis \cite{8}. Wang et al \cite{7} considered the $H_{\infty}$-controller design for NCSs with Markovian packet dropouts. Xie and Xie \cite{10} presented the necessary and sufficient condition for the mean-square stability of sampled-data networked linear systems with Markovian packet losses.\\
As said in \cite{14} and \cite{13}, the zero-input strategy is mainly for mathematical convenience rather than for performance considerations. Also, it must be pointed out that in most cases, especially in the transition phase, the hold-input strategy is better than the zero-input strategy. Hristu-Varsakelis \cite{12} analyzed the structure properties (e.g., observability and controllability) of the NCSs which a zero-order hold is included. Lu et al \cite{13} considered the NCSs with packet dropout using an improved switching hold compensation strategy in which the too old held signal is deleted. The stability and the optimal controller depending on the packet loss probability were obtained. Sun et al \cite{11} studied the $L_{2}$-gain of systems with input delays and controller temporary failure which is modeled as zero-order hold, and established the sufficient condition of exponential stability for the NCSs.\\
Most of the previous works only considered either delay or packet loss in NCSs. However, there are few works that focus on the simultaneous occurrence of these two uncertainties. In this paper, we consider the optimal control problem for NCSs simultaneously with input delay and Markovain dropout under a more complex, but more general and practical compensation strategy, i.e., hold-input strategy. To reduce the computational complexity caused by the hold-input strategy, using the state augmentation method, we firstly convert the system under hold-input strategy to the linear system with Markovian jump (MJLS), which is another important topic. For example, Li and Zhou \cite{07} and Li et al \cite{08} considered the indefinite stochastic optimal control problems for the MJLS over a finite time horizon and an infinite time horizon, respectively. Also, Costa et al  \cite{21} studied discrete-time Markovian jump linear systems and their applications, and Han et al  \cite{20} derived the optimal control for discrete-time Markovian jump linear system with control input delay.
On this basis,  the main obstacles to solve the stability of the system are the adaptability of controller and the temporal correlation caused by simultaneous input delay and Markovian jump.\\
In view of these, the key point in this paper is how to deal with the forward and backward stochastic difference equations (FBSDEs), which are derived by the stochastic maximum principle. Inspired by \cite{18} and \cite{19} in which the FBSDEs have made substantial progress in optimal LQ control problem for linear systems, the main results in this paper are derived and can be summarized as follows. First, the necessary and sufficient conditions for the solvability of optimal control problem over a finite horizon are presented by the coupled difference Riccati-type equations (CDREs). Second, the existence of the solution to the coupled algebraic Riccati-type equations (CAREs) is proved. Moreover, the optimal controller and optimal cost functional over an infinite horizon are derived. Finally, the necessary and sufficient conditions for the stabilization of the NCSs are established using the CAREs.\\
The rest of this article is structured as follows. Section 2 gives the problem statement. Section 3 solves the  optimal control problem over a finite horizon and the stabilization problems for infinite horizon case. An numerical example is presented to verify the obtained results in Section 4. A summary is presented in Section 5. Proofs for some results can be found in Appendix.\\
\emph{Notation }: \ ${\mathbb{R}}^n$ is the $n$-dimensional Euclidean space and $\mathbb{R}^{m\times n}$  the norm bounded linear space of all
$m\times n$ matrices. $Y'$ is  the transposition of $Y$ and  $Y\geq 0 (Y>0)$ means  that  $Y\in \mathbb{R}^{n\times n}$ is symmetric positive semi-definite  (positive definite). Let $(\Omega,\mathcal{F}, \mathcal{F}_{k}, \mathcal{P})$ be a complete probability space with the natural filtration  $\{\mathcal{F}_{k}\}_{k\geq0}$  generated by $\{\theta_{0},\cdots,\theta_{k}\}$. $\mathbb{E}[\cdot|\mathcal{F}_{k}]$ means  the conditional expectation with respect
to $\mathcal{F}_{k}$ and $\mathcal{F}_{-1}$ is understood as $\{\emptyset,\Omega\}$.

\section{Problem Statement and Preliminaries}
Consider the following discrete-time system:
\begin{eqnarray}
&&x_{k+1}=Ax_{k}+Bu^{a}_{k-d},\label{f1}\\
&&\left\{
\begin{array}{lll}
u^{a}_{k}=\theta_{k}u^{c}_{k}+(1-\theta_{k})u^{a}_{k-1}, \  k\geq 0,\\
u^{a}_{i}=u^{c}_{i}, \ i=-d,\cdots,-1,
\end{array}
\right.\label{f2}
\end{eqnarray}
where $x_{k}\in \mathbb{R}^{n}$ is the state. $u^{a}_{k}\in \mathbb{R}^{m}$ denotes the control input to the actuator and $u^{c}_{k}\in \mathbb{R}^{m}$ is the desired control input computed by the controller. The stochastic variable $\theta_{k}$ is packet dropout modeled as a two state Markov chain $\theta_{k}\in\{0,1\}$ with transition probability  $\xi_{ij}=\mbox{P}(\theta_{k+1}=j|\theta_{k}=i)(i,j=0,1)$ between the controller and the actuator: Take $u^{a}_{k}=u^{c}_{k}$, if the packet is correctly delivered; otherwise, take $u^{a}_{k}=u^{a}_{k-1}$, if the packet is lost. The initial values are $x_{0}, u^{c}_{-d}, \cdots, u^{c}_{-1}$. In addition, $A\in \mathbb{R}^{n\times n}$, $B\in \mathbb{R}^{n\times m}$ are constant matrices.\\
By state augmentation, we have
\begin{eqnarray}
\left[
\begin{array}{cc}
x_{k+1}\\
u^{a}_{k-d}
\end{array}\right]\hspace{-1mm}=\hspace{-1mm}\left[
\begin{array}{cc}
A&(1-\theta_{k})B\\
0&(1-\theta_{k})I
\end{array}\right]\left[
\begin{array}{cc}
x_{k}\\
u^{a}_{k-d-1}
\end{array}\right]\hspace{-1mm}+\hspace{-1mm}
\left[
\begin{array}{cc}
\theta_{k}B\\
\theta_{k}I
\end{array}\right]u^{c}_{k-d}.\label{f4}
\end{eqnarray}
Let $z_{k+1}=\left[
\begin{array}{cc}
x_{k+1}\\
u^{a}_{k-d}
\end{array}\right]$, $\bar{A}_{\theta_{k}}=\left[
\begin{array}{cc}
A&(1-\theta_{k})B\\
0&(1-\theta_{k})I
\end{array}\right]$, $\bar{B}_{\theta_{k}}=\left[
\begin{array}{cc}
\theta_{k}B\\
\theta_{k}I
\end{array}\right]$, then (\ref{f4}) can be rewritten as
\begin{eqnarray}
z_{k+1}=\bar{A}_{\theta_{k}}z_{k}+\bar{B}_{\theta_{k}}u^{c}_{k-d}. \label{f5}
\end{eqnarray}\\
{\bf Remark 1 \label{r1} } Obviously, systems (\ref{f1}) and (\ref{f2}) are equivalent to the augmented systems (\ref{f5}), i.e., the linear system with simultaneous input delay and Markovian jump. \\
Define the following cost functional over an infinite horizon as
\begin{eqnarray}
J=\mathbb{E}\Bigg[\sum^{\infty}_{k=0}\Big(z_{k}'Qz_{k}
+(u^{c}_{k-d})'Ru^{c}_{k-d}\Big)\Bigg], \label{fi01}
\end{eqnarray}
where weighting matrices $Q\in \mathbb{R}^{(n+m)\times (n+m)}, R\in \mathbb{R}^{m\times m}$.\\
{\bf Problem 1 \label{p3}} Find a $\mathcal{F}_{k-1}$-measurable controller $u^{c}_{k}$ to stabilize system (\ref{f1})-(\ref{f2}) and minimize cost functional (\ref{fi01}).\\
{\bf Remark 2\label{r01} }  Most of works in the literatures only considered either delay or packet loss in NCSs. However, there are few works that focus on the simultaneous occurrence of these two uncertainties. In this paper, the optimal control problem for NCSs including both input delay and Markovain dropout under a more general yet practical compensation strategy (for example, hold-input strategy) will be considered. Due to the  adaptability of controller and the temporal correlation caused by simultaneous input delay and Markovian jump, there exists challenging to solve the stability of the system.

\section{Main Results}
For discussion, this section will follow two steps. The LQ optimal control problem over a finite horizon will be first considered by solving the FBSDEs which is very vital in this paper. On this basis, Problem 1 will be resolved.
\subsection{LQ Optimal Control in Finite Horizon Case}
Consider the following cost functional:
\begin{eqnarray}
J_{N}&=&\mathbb{E}\Bigg[\sum^{N}_{k=0}\Big(z_{k}'Qz_{k}
+(u^{c}_{k-d})'Ru^{c}_{k-d}\Big)\nonumber\\
&&+z_{N+1}'\bar{P}_{N+1}z_{N+1}\Bigg], \label{f3}
\end{eqnarray}
where weighting matrices $Q\in \mathbb{R}^{(n+m)\times (n+m)}, R\in \mathbb{R}^{m\times m}$ and terminal value $\bar{P}_{N+1}$ are positive semi-definite.\\
{\bf Problem 2 } Find a $\mathcal{F}_{k-d-1}$-measurable controller $u^{c}_{k-d}$ such that cost functional (\ref{f3}) is minimized subject to (\ref{f5}).\\
Applying the maximum principle to Problem 2, the following FBSDEs are obtained
\begin{eqnarray}
\left\{
\begin{array}{lll}
0=\mathbb{E}_{k-d-1}\big[\bar{B}_{\theta_{k}}'\lambda_{k}\big]
+Ru^{c}_{k-d},\\
\lambda_{k-1}=Qz_{k}
+\mathbb{E}_{k-1}\big[\bar{A}_{\theta_{k}}'\lambda_{k}\big],\\
\lambda_{N}=\bar{P}_{N+1}z_{N+1},\\
z_{k+1}=\bar{A}_{\theta_{k}}z_{k}+\bar{B}_{\theta_{k}}u^{c}_{k-d}.
\end{array}
\right.\label{f07}
\end{eqnarray}
{\bf Remark 3 }    Due to the Markovian jump and input delay which gives rise to the fundamental difficulty about the adaptability of controller and the temporal correlation, our problem is important and challenging compared with \cite{18} and \cite{19}. The key technique in this paper is to tackle FBSDEs (\ref{f07}). Further, we can solve our problem via FBSDEs (\ref{f07}). \\
For any $d \leq k\leq N$, define the following recursive sequence,
\begin{eqnarray}
\bar{P}_{\theta_{k-1}}(k)&=&Q+\mathbb{E}_{k-1}\big[(\bar{A}_{\theta_{k}})'
\bar{P}_{\theta_{k}}(k+1)\bar{A}_{\theta_{k}}\nonumber\\
&&-(M^{0}_{\theta_{k-1}})'
\Gamma^{-1}_{\theta_{k-1}}M^{0}_{\theta_{k-1}}\big],\label{f8}
\end{eqnarray}
where
\begin{eqnarray}
\Gamma_{\theta_{k-d-1}}\hspace{-1mm}&=&\hspace{-1mm}R\hspace{-1mm}+\hspace{-1mm}
\mathbb{E}_{k-d-1}\Bigg[(\bar{B}_{\theta_{k}})'
\bar{P}_{\theta_{k}}(k+1)\bar{B}_{\theta_{k}}\hspace{-1mm}-\hspace{-1mm}\sum^{d-1}_{i=0}(M^{i+1}_{\theta_{k-d+i}})'\nonumber\\
\hspace{-1mm}&&\times\hspace{-1mm} \Gamma^{-1}_{\theta_{k-d+i}}M^{i+1}_{\theta_{k-d+i}}\Bigg],\label{f9}\\
M^{0}_{\theta_{k-d-1}}\hspace{-1mm}&=&\hspace{-1mm}\mathbb{E}_{k-d-1}
\Bigg[(\tilde{S}^{1}_{\theta_{k-1}})'
\prod^{d}_{j=1}\bar{A}_{\theta_{k-j}}-\sum^{d-1}_{i=0}\bigg((M^{i+1}_{\theta_{k-d+i}})'\nonumber\\
\hspace{-1mm}&&\times\hspace{-1mm}
\Gamma^{-1}_{\theta_{k-d+i}}M^{0}_{\theta_{k-d+i}}\prod^{i}_{s=0}\bar{A}_{\theta_{k-d+s}}\bigg)\Bigg],
\label{f10}\\
M^{i}_{\theta_{k-d-1}}\hspace{-1mm}&=&\hspace{-1mm}\mathbb{E}_{k-d-1}
\Bigg[(\tilde{S}^{1}_{\theta_{k-1}})'
\prod^{i-1}_{j=1}\bar{A}_{\theta_{k-j}}\bar{B}_{\theta_{k-i}}
\hspace{-1mm}-\hspace{-1mm}\sum^{d-1}_{s=0}(M^{s+1}_{\theta_{k-d+s}})'
\nonumber\\
\hspace{-1mm}&&\times\hspace{-1mm} \Gamma^{-1}_{\theta_{k-d+s}}M^{i+s+1}_{\theta_{k-d+s}}\Bigg],\ \ i=1,\cdots,d,
\label{f1100}\\
M^{i}_{\theta_{k-d-1}}\hspace{-1mm}&=&\hspace{-1mm}M^{0}_{\theta_{k-d-1}}\mathbb{E}_{k-d-1}\Bigg[\prod^{i-d-1}_{j=1}\bar{A}_{\theta_{k-d-j}}
\bar{B}_{\theta_{k-i}}\Bigg], \label{f0011}\\
 \hspace{-1mm}&& \hspace{-1mm}   \ \ \ \ \ \ \ \ \ \ i\geq d+1, \nonumber\\
 M^{i}_{\theta_{N-s-1}}\hspace{-1mm}&=&\hspace{-1mm}0, \ \ \  i\geq0, \ \ \ s\leq d-1,\label{f011}\\
\tilde{S}^{1}_{\theta_{k-1}}\hspace{-1mm}&=&\hspace{-1mm}\mathbb{E}_{k-1}\big[(\bar{A}_{\theta_{k}})'
\bar{P}_{\theta_{k}}(k+1)\bar{B}_{\theta_{k}}\big],\label{f0012}\\
\tilde{S}^{j}_{\theta_{k-1}}\hspace{-1mm}&=&\hspace{-1mm}\mathbb{E}_{k-1}\big[(\tilde{S}^{j-1}_{\theta_{k}})'
\bar{A}_{\theta_{k}}\big],\label{f12}
\end{eqnarray}
with terminal value $\bar{P}_{\theta_{N}}(N+1)=\bar{P}_{N+1}$. And equations (\ref{f8})-(\ref{f12}) are termed the CDREs.\\
{\bf Remark 4} Using the ineffectiveness theory of Markovian, the conditional expectation of the sequence $f_{\theta_{k}}$ can be expressed as
\begin{eqnarray}
\hspace{-1mm}&&\hspace{-1mm}\mathbb{E}_{k-j}[f_{\theta_{k}}]\hspace{-1mm}=\hspace{-1mm}\sum^{1}_{\theta_{k-j+1}=0}\xi_{\theta_{k-j+1}\theta_{k-j}}\Bigg\{
\sum^{1}_{\theta_{k-j+2}=0}\xi_{\theta_{k-j+2}\theta_{k-j+1}}\cdots\nonumber\\
\hspace{-1mm}&&\times\hspace{-1mm}
\Bigg[\sum^{1}_{\theta_{k-1}=0}\xi_{\theta_{k-1}\theta_{k-2}}
(\sum^{1}_{\theta_{k}=0}\xi_{\theta_{k}\theta_{k-1}})f_{\theta_{k}}\Bigg]\Bigg\}.\label{f13}
\end{eqnarray}\\
{\bf Remark 5} For convenience, set  $\Gamma_{\theta_{k-d-1}}(k)\triangleq\Gamma_{\theta_{k-d-1}},
\bar{P}_{\theta_{k}}(k+1)\triangleq\bar{P}_{\theta_{k}},M^{i}_{\theta_{k-d-1}}(k)\triangleq M^{i}_{\theta_{k-d-1}}$, $\tilde{S}^{1}_{\theta_{k-1}}(k)\triangleq \tilde{S}^{1}_{\theta_{k-1}}$. \\
{\bf Lemma 1 } From CDREs (\ref{f8})-(\ref{f12}), the following relationships can be given
\begin{eqnarray}
M^{0}_{\theta_{k-d-1}}\hspace{-1mm}&=&\hspace{-1mm}\mathbb{E}_{k-d-1}\Bigg[(F^{1}_{\theta_{k-1}})'\prod^{d}_{j=1}
\bar{A}_{\theta_{k-j}}\Bigg],\label{e02}\\
M^{i}_{\theta_{k-d-1}}\hspace{-1mm}&=&\hspace{-1mm}\mathbb{E}_{k-d-1}\big[(F^{i}_{\theta_{k-i}})'\bar{B}_{\theta_{k-i}}\big], \ \ \ i\geq1,\label{e0}\\
F^{i}_{\theta_{k-1}}\hspace{-1mm}&=&\hspace{-1mm}\mathbb{E}_{k-1}\Bigg[\prod^{i-2}_{j=0}(\bar{A}_{\theta_{k+j}})'
F^{1}_{\theta_{k+i-2}}\Bigg],\ \ \ i\geq 2, \label{e01}
\end{eqnarray}
with
\begin{eqnarray}
F^{1}_{\theta_{k-1}}\hspace{-1mm}&=&\hspace{-1mm}\tilde{S}^{1}_{\theta_{k-1}}-\sum^{d}_{i=0}F^{i+2}_{\theta_{k-1}}
\Gamma^{-1}_{\theta_{k-d+i}}M^{i+1}_{\theta_{k-d+i}},\label{f17}\\
F^{i}_{\theta_{k-1}}\hspace{-1mm}&=&\hspace{-1mm}\mathbb{E}_{k-1}[(\bar{A}_{\theta_{k}})'F^{i-1}_{\theta_{k}}],\ \label{f18}\\
F^{N-k+1}_{\theta_{k-1}}&=&(\tilde{S}^{N-k+1}_{\theta_{k-1}})'.\label{f018}
\end{eqnarray}
Proof. From (\ref{f9})-(\ref{f12}) and via mathematical induction, (\ref{e02})-(\ref{e01}) can be simply calculated, here, we omit the proof. \\
Based on the preliminaries, the results of Problem 2 can be obtained in this section.\\
{\bf Theorem 1}  There exists a unique solution to Problem 2 if and only if  $\Gamma_{\theta_{k-d-1}}$ in (\ref{f8}) is positive definite. In this case, the optimal controller can be given as
\begin{eqnarray}
u^{c}_{k-d}\hspace{-1mm}=\hspace{-1mm}-\Gamma^{-1}_{\theta_{k\hspace{-0.5mm}-\hspace{-0.5mm}
d\hspace{-0.5mm}-\hspace{-0.5mm}1}}\Bigg(M^{0}_{\theta_{k-d-1}}z_{k-d}
\hspace{-1mm}+\hspace{-1mm}\sum^{d}_{i=1}M^{i}_{\theta_{k\hspace{-0.5mm}-\hspace{-0.5mm}d
\hspace{-0.5mm}-\hspace{-0.5mm}1}}u^{c}_{k-d-i}\Bigg), \label{f14}
\end{eqnarray}
and the optimal cost functional is
\begin{eqnarray}
J_{N}\hspace{-1mm}&=&\hspace{-1mm}\mathbb{E}\Bigg\{\sum_{k=0}^{d-1}\big[z_{k}'Qz_{k}
+(u^{c}_{k-d})'Ru^{c}_{k-d}\big]
+z_{d}'\bar{P}_{\theta_{d-1}}z_{d}\nonumber\\
&&\hspace{-1mm}-\hspace{-1mm}z_{d}'\sum^{d-1}_{s=0}\left(F^{s+1}_{\theta_{d-1}}
\Gamma^{-1}_{\theta_{s-1}}M^{0}_{\theta_{s-1}}z_{s}\right)\nonumber\\
&&\hspace{-1mm}-\hspace{-1mm}z_{d}'\sum^{d}_{s=0}\left(F^{s+1}_{\theta_{d-1}}
\Gamma^{-1}_{\theta_{s-1}}\sum^{d}_{i=s+1}M^{i}_{\theta_{s-1}}u^{c}_{s-i}\right)\Bigg\}. \label{f15}
\end{eqnarray}
Moreover, the solution of the FBSDEs can be given
\begin{eqnarray}
\lambda_{k-1}\hspace{-1mm}&=&\hspace{-1mm}\bar{P}_{\theta_{k-1}}z_{k}\hspace{-1mm}-\hspace{-1mm}
\sum^{d-1}_{s=0}\left(F^{s+1}_{\theta_{k-1}}
\Gamma^{-1}_{\theta_{k\hspace{-0.5mm}-\hspace{-0.5mm}d\hspace{-0.5mm}-\hspace{-0.5mm}1\hspace{-0.5mm}
+\hspace{-0.5mm}s}}M^{0}_{\theta_{k\hspace{-0.5mm}-\hspace{-0.5mm}d\hspace{-0.5mm}-\hspace{-0.5mm}1\hspace{-0.5mm}
+\hspace{-0.5mm}s}}z_{k-d+s}\right)\nonumber\\
&&\hspace{-1mm}-\hspace{-1mm}\sum^{d-1}_{s=0}\left(F^{s+1}_{\theta_{k-1}}
\Gamma^{-1}_{\theta_{k\hspace{-0.5mm}-\hspace{-0.5mm}d\hspace{-0.5mm}-\hspace{-0.5mm}1\hspace{-0.5mm}
+\hspace{-0.5mm}s}}\sum^{d}_{i=s+1}M^{i}_{\theta_{k\hspace{-0.5mm}-\hspace{-0.5mm}d\hspace{-0.5mm}-\hspace{-0.5mm}1\hspace{-0.5mm}
+\hspace{-0.5mm}s}}u^{c}_{k-d-i+s}\right).\label{f16}
\end{eqnarray}
Proof. See Appendix \ref{ap2}.\\
{\bf Remark 6} For convenience, the notation has been written for short, i.e., $F^{i}_{\theta_{k-1}}(k)\triangleq F^{i}_{\theta_{k-1}}$.\\
{\bf Remark 7}
 For the delay-free case, i.e., d=0 in systems (\ref{f1})-(\ref{f2}),  then the result of Theorem 1 can be rewritten as follows. The optimal controller (\ref{f14}) and the solution of the FBSDEs (\ref{f16}) can be re-expressed as
 \begin{eqnarray}
u^{c}_{k}&=&-\Gamma^{-1}_{\theta_{k-1}}M^{0}_{\theta_{k-1}}z_{k}, \label{b1}\\
\lambda_{k-1}&=&\bar{P}_{\theta_{k-1}}z_{k}, \label{b2}
\end{eqnarray}
where
\begin{eqnarray}
\bar{P}_{\theta_{k-1}}(k)&=&Q+\mathbb{E}_{k-1}\big[(\bar{A}_{\theta_{k}})'
\bar{P}_{\theta_{k}}(k+1)\bar{A}_{\theta_{k}}\nonumber\\
&&-(M^{0}_{\theta_{k-1}})'
\Gamma^{-1}_{\theta_{k-1}}M^{0}_{\theta_{k-1}}\big],\label{b3}\\
\Gamma_{\theta_{k-1}}&=&R+\mathbb{E}_{k-1}[(\bar{B}_{\theta_{k}})'
\bar{P}_{\theta_{k}}(k+1)\bar{B}_{\theta_{k}}],\label{b9}\\
M^{0}_{\theta_{k-1}}&=&\mathbb{E}_{k-1}[(\bar{B}_{\theta_{k}})'
\bar{P}_{\theta_{k}}(k+1)\bar{B}_{\theta_{k}}],\label{b10}
\end{eqnarray}
which are parallel to the results of standard case for MJLSs \cite{0012}.
\subsection{Optimal Control in Infinite Horizon Case}
In the sequence, the following assumptions will be made.\\
{\bf Assumption 1 \label{a1}}  The state weighting matrix $Q$ is positive semi-definite and the control weighting matrix $R$ is strictly positive definite.\\
{\bf Assumption 2 \label{a2}} ($\bar{A}, Q^{\frac{1}{2}}$) is exactly observable, where $\bar{A}=(\bar{A}_{0}, \bar{A}_{1})$.\\
Define the following CAREs for $l_{j}=0,1, 0\leq j\leq d$,
\begin{eqnarray}
\bar{P}_{l_{d}}\hspace{-1mm}=\hspace{-1mm}Q+\mathbb{E}_{l_{d}}\big[(\bar{A}_{l_{d+1}})'
\bar{P}_{l_{d+1}}\bar{A}_{l_{d+1}}-(M^{0}_{l_{d}})'
\Gamma^{-1}_{l_{d}}M^{0}_{l_{d}}\big],\label{f17}
\end{eqnarray}
where
\begin{eqnarray}
\Gamma_{l_{0}}\hspace{-1mm}&=&\hspace{-1mm}R\hspace{-1mm}+\hspace{-1mm}\mathbb{E}_{l_{0}}
\Bigg[(\bar{B}_{l_{d\hspace{-0.51mm}+\hspace{-0.51mm}1}})'
\bar{P}_{l_{d\hspace{-0.51mm}+\hspace{-0.51mm}1}}\bar{B}_{l_{d\hspace{-0.51mm}+\hspace{-0.51mm}1}}
\hspace{-1.5mm}-\hspace{-1.5mm}\sum^{d}_{i=0}
(M^{i\hspace{-0.51mm}+\hspace{-0.51mm}1}_{l_{i\hspace{-0.51mm}+\hspace{-0.51mm}1}})'
\Gamma^{-1}_{l_{i\hspace{-0.51mm}+\hspace{-0.51mm}1}}
M^{i\hspace{-0.51mm}+\hspace{-0.51mm}1}_{l_{i\hspace{-0.51mm}+\hspace{-0.51mm}1}}\Bigg],\label{f18}\\
M^{0}_{l_{0}}\hspace{-1mm}&=&\hspace{-1mm}\mathbb{E}_{l_{0}}\Bigg[(\tilde{S}^{1}_{l_{d}})'
\prod^{d}_{j=1}\bar{A}_{l_{d\hspace{-0.51mm}+\hspace{-0.51mm}1\hspace{-0.51mm}-\hspace{-0.51mm}j}}
\hspace{-1.5mm}-\hspace{-1.5mm}\sum^{d}_{i=0}
(M^{i\hspace{-0.51mm}+\hspace{-0.51mm}1}_{l_{i\hspace{-0.51mm}+\hspace{-0.51mm}1}})'
\Gamma^{-1}_{l_{i\hspace{-0.51mm}+\hspace{-0.51mm}1}}M^{0}_{l_{i\hspace{-0.51mm}+\hspace{-0.51mm}1}}
\prod^{i}_{s=0}\bar{A}_{l_{s\hspace{-0.51mm}+\hspace{-0.51mm}1}})\Bigg],\label{f19}\\
M^{i}_{l_{0}}\hspace{-1mm}&=&\hspace{-1mm}\mathbb{E}_{l_{0}}\Bigg[(\tilde{S}^{1}_{l_{d}})'
\prod^{i-1}_{j=1}\bar{A}_{l_{d\hspace{-0.51mm}+\hspace{-0.51mm}1\hspace{-0.51mm}-\hspace{-0.51mm}j}}
\bar{B}_{l_{d\hspace{-0.51mm}+\hspace{-0.51mm}1\hspace{-0.51mm}-\hspace{-0.51mm}i}}
\hspace{-1.5mm}-\hspace{-1.5mm}\sum^{d}_{s=0}
(M^{s\hspace{-0.51mm}+\hspace{-0.51mm}1}_{l_{s\hspace{-0.51mm}+\hspace{-0.51mm}1}})'
\Gamma^{-1}_{l_{s\hspace{-0.51mm}+\hspace{-0.51mm}1}}
M^{i\hspace{-0.51mm}+\hspace{-0.51mm}s\hspace{-0.51mm}+\hspace{-0.51mm}1}_{l_{s\hspace{-0.51mm}+\hspace{-0.51mm}1}}\Bigg],\label{f20}\\
&&  \ \ \ \ \ i=1,\cdots,d,\nonumber\\ M^{s}_{l_{1}}\hspace{-1mm}&=&\hspace{-1mm}M^{0}_{l_{1}}\mathbb{E}_{l_{1}}\Bigg[\prod^{s-d-1}_{j=1}\bar{A}_{l_{j}}
\bar{B}_{l_{s-d}}\Bigg], \ \ \ s\geq d+1 \label{f0020}\\
\tilde{S}^{1}_{l_{d}}\hspace{-1mm}&=&\hspace{-1mm}\mathbb{E}_{l_{d}}\big[(\bar{A}_{l_{d+1}})'
\bar{P}_{l_{d+1}}\bar{B}_{l_{d+1}}\big],\ \ \label{f21}\\
\tilde{S}^{j}_{l_{d}}\hspace{-1mm}&=&\hspace{-1mm}\mathbb{E}_{l_{d}}\big[(\bar{A}_{l_{d+1}})'\tilde{S}^{j-1}_{l_{d+1}}\big].\label{f021}
\end{eqnarray}
The main result will be presented next.\\
{\bf Theorem 2 }
Under  Assumption 1 and 2, the system (\ref{f1}) is stabilizable in the mean-square sense if and only if CAREs (\ref{f17})-(\ref{f21}) have a solution such that
\begin{eqnarray}
\bar{P}_{l_{d}}
-\sum_{s=0}^{d-1}\big[(F^{s+1}_{l_{d}})'
\Gamma^{-1}_{l_{s-1}}
F^{s+1}_{l_{d}}\big]>0,\label{f22}
\end{eqnarray}
in which
\begin{eqnarray}
F^{1}_{l_{d}}&=&\tilde{S}^{1}_{l_{d}}-\sum^{d}_{i=0}F^{i+2}_{l_{d}}
\Gamma^{-1}_{l_{i+1}}M^{i+1}_{l_{i+1}},\label{f23}\\
F^{i}_{l_{d}}&=&\mathbb{E}_{l_{d}}\big[(\bar{A}_{l_{d+1}})'F^{i-1}_{l_{d+1}}\big],\ \
F^{d}_{l_{d}}=(\tilde{S}^{d}_{l_{d}})',\label{f24}
\end{eqnarray}
$ l_{i}\in \{0,1\}$, $i=0,1,\cdots,d+1$.
Moreover, for $k\geq d$ the optimal controller can be given as
\begin{eqnarray}
u^{c}_{k-d}&=&-\Gamma^{-1}_{l_{0}}\Bigg(M^{0}_{l_{0}}z_{k-d}
+\sum^{d}_{i=1}M^{i}_{l_{0}}u^{c}_{k-d-i}\Bigg). \label{f25}
\end{eqnarray}
 The corresponding cost index is presented by
\begin{eqnarray}
J^{\ast}\hspace{-1mm}&=&\hspace{-1mm}\mathbb{E}\Bigg[z_{0}'\bar{P}_{l_{d}}z_{0}
\hspace{-1mm}+\hspace{-1mm}\sum_{k=0}^{d-1}\Big(u^{c}_{k\hspace{-0.51mm}-\hspace{-0.51mm}d}
\hspace{-1mm}+\hspace{-1mm}\Gamma^{-1}_{i}M^{0}_{i}z_{k\hspace{-0.51mm}-\hspace{-0.51mm}d}
\hspace{-1mm}+\hspace{-1mm}\Gamma^{-1}_{i}\nonumber\\
&&\times \sum^{d}_{s=1}(M^{s}_{i}u^{c}_{k\hspace{-0.51mm}
-\hspace{-0.51mm}d\hspace{-0.51mm}-\hspace{-0.51mm}s})\Big)'\Gamma_{i}\Big(u^{c}_{k\hspace{-0.51mm}-\hspace{-0.51mm}d}\hspace{-1mm}
+\hspace{-1mm}\Gamma^{-1}_{i}M^{0}_{i}z_{k\hspace{-0.51mm}-\hspace{-0.51mm}d}\nonumber\\
&&
\hspace{-1mm}+\hspace{-0.5mm}\Gamma^{-1}_{i}\sum^{d}_{s=1}(M^{s}_{i}u^{c}_{k\hspace{-0.51mm}
-\hspace{-0.51mm}d\hspace{-0.51mm}-\hspace{-0.51mm}s})\Big)\Bigg],\label{f26}
\end{eqnarray}
where $\bar{P}_{s}, \Gamma_{s}, M^{i}_{s}, s=0,1$ satisfy CAREs (\ref{f17})-(\ref{f21}).\\
Proof. See Appendix \ref{ap3}.\\
{\bf Remark 8} The optimal control for the MJLS without delay has been well studied in the literature, and using the state augmentation method, we resolve the presented problem in this paper. However, it will bring large amount of calculation especially for the high dimension system or the large delay.\\
{\bf Remark 9} Compared with the previous works only considered either delay or packet loss or under zero-input strategy in NCSs (\cite{Imer: 06}, \cite{13} and so on), the necessary and sufficient conditions for the stabilization of the NCSs including both input delay and Markovain dropout under hold-input strategy are established. To the best of our knowledge, the above necessary and sufficient conditions are firstly presented. \\
\section{Numerical Example}
Consider system (\ref{f1}) with $A=1, B=15, d=1$, and the initial values $x_{0}=10, u^{c}_{-1}=1$, let the transition probability $\xi_{00}=0.9, \xi_{11}=0.7$ and cost functional (\ref{f3}) with $Q=\left[
  \begin{array}{cc}
1&0\\
0&1
  \end{array}
\right], R=10 $. Therefore,
$\bar{A}_{0}=\left[
\begin{array}{cc}
1&15\\
0&1
\end{array}\right], \bar{A}_{1}=\left[
\begin{array}{cc}
1&0\\
0&0
\end{array}\right], \bar{B}_{0}=\left[
\begin{array}{cc}
0\\
0
\end{array}\right], \bar{B}_{1}=\left[
\begin{array}{cc}
15\\
1
\end{array}\right]$.
In this case, a sample path of the Markov chain $\theta_{k}$ is shown in Figure 1.\\
In view of (\ref{f17})-(\ref{f21}), the following results can be obtained:
 \begin{eqnarray*}
P_{0}&=&\left[
  \begin{array}{cc}
3.7383&251.49\\
251.49&69116.31
  \end{array}
\right], P_{1}=\left[
  \begin{array}{cc}
3.8049&88.0818\\
88.0818&23233.77
  \end{array}
\right], \\
\Gamma_{0}&=&3006.02, \Gamma_{1}=11877.83,\\
 M^{0}_{0}&=&\left[
  \begin{array}{cc}
55.01&1383.48
  \end{array}
\right],
M^{0}_{1}=\left[
  \begin{array}{cc}
107.89&461.16
  \end{array}
\right], \\
M^{1}_{0}&=&-1834.41, M^{1}_{1}=-5058.20.
\end{eqnarray*}
Hence, (\ref{f22}) can be calculated
\begin{eqnarray}
\bar{P}_{0}
-(F^{1}_{0})'
\Gamma^{-1}_{i}
F^{1}_{0}&=&\left[
  \begin{array}{cc}
4  &  251 \\
    251  &  69116
  \end{array}
\right]>0,i=0,1,\\
 \bar{P}_{1}
-(F^{1}_{1})'
\Gamma^{-1}_{0}
F^{1}_{1}&=&\left[
  \begin{array}{cc}
0.37  &  88 \\
    88  &  23234
  \end{array}
\right]>0,\\
\bar{P}_{1}
-(F^{1}_{1})'
\Gamma^{-1}_{1}
F^{1}_{1}&=&\left[
  \begin{array}{cc}
3  &  88 \\
    88  &  23234
  \end{array}
\right]>0.
\end{eqnarray}
According to Theorem 2, the optimal controller can be expressed as
\begin{eqnarray*}
u^{c}_{k}&=&-\left[
  \begin{array}{cc}
0.0183&0.4602
  \end{array}
\right]z_{k}+0.6102u^{c}_{k-1}, l_{0}=0,\\
u^{c}_{k}&=&-\left[
  \begin{array}{cc}
0.0091&0.0388
  \end{array}
\right]z_{k}+0.4259u^{c}_{k-1}, l_{0}=1.
\end{eqnarray*}
 A simulation result of optimal controller $u^{c}_{k}$ and $u^{a}_{k}$ is shown in Figure 2 (a).\\
 \begin{figure}[htbp]
  \begin{center}
  \includegraphics[width=0.42\textwidth]{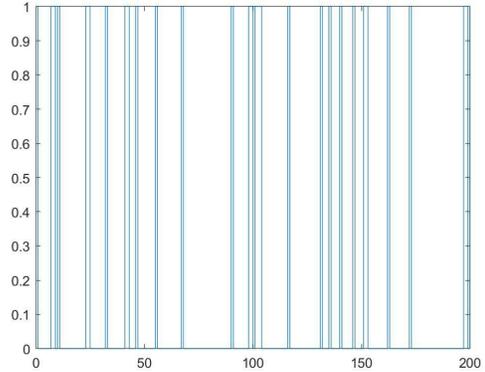}
  \caption{A sample path with q=0.9 and p=0.7} \label{fig:digit}
  \end{center}
\end{figure}
\begin{figure}[htbp]
  \begin{center}
  \includegraphics[width=0.42\textwidth]{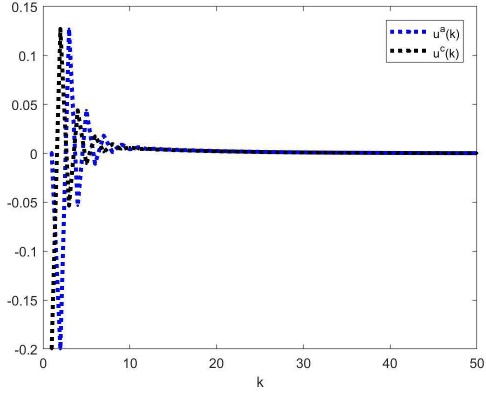}
  \caption{Optimal control.} \label{fig:digit}
  \end{center}
\end{figure}
\begin{figure}[htbp]
  \begin{center}
  \includegraphics[width=0.42\textwidth]{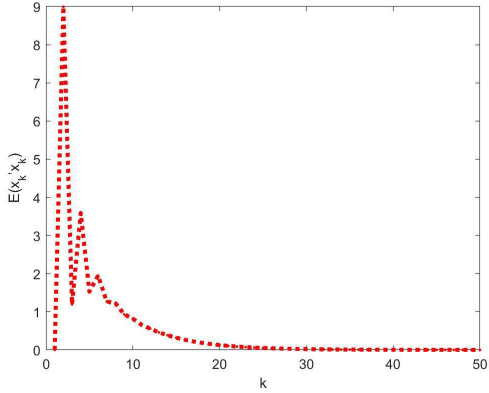}
  \caption{The trajectory of $E[x_{k}'x_{k}]$.} \label{fig:digit}
  \end{center}
\end{figure}
From Figure 2 (b), it is easy to see that system (\ref{f1})-(\ref{f2}) is stabilizable in the mean-square sense.
\section{Conclusion}
In this paper, the optimal LQ control problem for NCSs simultaneously with input delay and Markovian dropout is discussed. Compared with the results in the literature, we mainly consider the hold-input strategy, which is much more computationally complicated than zero-input strategy. Necessary and sufficient conditions for the solvability of optimal control problem over a finite horizon are presented by the CDREs. Moreover, the NCS is mean-square stability if and only if the CAREs have a particular solution. The key technique in this paper is to tackle the FBSDEs, which are more difficult to be dealt with, due to the adaptability of controller and the temporal correlation caused by simultaneous input delay and Markovian jump.

\appendices
\section{ Proof of Theorem 1 \label{ap2}}
Proof. ($\Longrightarrow$) $\Gamma_{\theta_{k-d-1}}>0$ will be proved by mathematical induction. Denote
\begin{eqnarray}
\tilde{J}_{k}\hspace{-1mm}&=&\hspace{-1mm}\mathbb{E}\Bigg\{\sum_{i=k}^N\big(z_{i}'Qz_{i}
\hspace{-1mm}+\hspace{-1mm}(u^{c}_{i\hspace{-0.51mm}-\hspace{-0.51mm}d})'R
u^{c}_{i\hspace{-0.51mm}-\hspace{-0.51mm}d}\big)\hspace{-1mm}+\hspace{-1mm}z_{N+1}'
\bar{P}_{N\hspace{-0.51mm}+\hspace{-0.51mm}1}z_{N\hspace{-0.51mm}+\hspace{-0.51mm}1}\Bigg\}. \label{f100}
\end{eqnarray}
Let $k=N$ in (\ref{f100}) with $z_{N}=0$, we have
\begin{eqnarray}
\tilde{J}_{N}&=&\mathbb{E}\big[(u^{c}_{N-d})'(R+\bar{B}_{\theta_{N}}'
\bar{P}_{N+1}\bar{B}_{\theta_{N}})u^{c}_{N-d}\big]\nonumber\\
&=&\mathbb{E}\big[(u^{c}_{N-d})'\Gamma_{\theta_{N-d-1}}u^{c}_{N-d}\big].\label{f0100}
\end{eqnarray}
Due to the uniqueness of the solution to Problem 2 and the arbitrariness of nonzero $u^{c}_{N-d}$, we obtain $\Gamma_{\theta_{N-d-1}}>0$. In this case, it follows from (\ref{f07}) that $u^{c}_{N-d}$ and $\lambda_{N-1}$ can be calculated as follows
\begin{eqnarray}
0\hspace{-1mm}&=&\hspace{-1mm}\mathbb{E}_{N\hspace{-0.51mm}-\hspace{-0.51mm}d
\hspace{-0.51mm}-\hspace{-0.51mm}1}\big[\bar{B}_{\theta_{N}}'\bar{P}_{N\hspace{-0.51mm}
+\hspace{-0.51mm}1}
(\bar{A}_{\theta_{N}}z_{N}+\bar{B}_{\theta_{N}}u^{c}_{N\hspace{-0.51mm}-\hspace{-0.51mm}d})
\hspace{-1mm}+\hspace{-1mm}
Ru^{c}_{N\hspace{-0.51mm}-\hspace{-0.51mm}d})\big]\nonumber\\
\hspace{-1mm}&=&\hspace{-1mm}\mathbb{E}_{N\hspace{-0.51mm}-\hspace{-0.51mm}d\hspace{-0.51mm}
-\hspace{-0.51mm}1}\big[\bar{B}_{\theta_{N}}'\bar{P}_{N\hspace{-0.51mm}+\hspace{-0.51mm}1}
\bar{A}_{\theta_{N}}z_{N}\hspace{-1mm}+\hspace{-1mm}(R+\bar{B}_{\theta_{N}}'\bar{P}_{N\hspace{-0.51mm}
+\hspace{-0.51mm}1}
\bar{B}_{\theta_{N}})u^{c}_{N\hspace{-0.51mm}-\hspace{-0.51mm}d}\big]\nonumber\\
\hspace{-1mm}&=&\hspace{-1mm}\mathbb{E}_{N\hspace{-0.51mm}-\hspace{-0.51mm}d
\hspace{-0.51mm}-\hspace{-0.51mm}1}\Bigg[(\tilde{S}^{1}_{\theta_{N\hspace{-0.51mm}-\hspace{-0.51mm}1}})'
\prod^{d}_{j=1}\bar{A}_{\theta_{N\hspace{-0.51mm}-\hspace{-0.51mm}j}}\Bigg]z_{N\hspace{-0.51mm}
-\hspace{-0.51mm}d}
\hspace{-1mm}+\hspace{-1mm}\sum^{d}_{i=1}\mathbb{E}_{N\hspace{-0.51mm}-\hspace{-0.51mm}
d\hspace{-0.51mm}-\hspace{-0.51mm}1}\Bigg[(\tilde{S}^{1}_{\theta_{N\hspace{-0.51mm}
-\hspace{-0.51mm}1}})'\nonumber\\
\hspace{-1mm}&&\times\hspace{-1mm}
\prod^{i-1}_{j=1}\bar{A}_{\theta_{N\hspace{-0.51mm}-\hspace{-0.51mm}j}}
\bar{B}_{\theta_{N\hspace{-0.51mm}-\hspace{-0.51mm}i}}\Bigg]u^{c}_{N\hspace{-0.51mm}
-\hspace{-0.51mm}d\hspace{-0.51mm}-\hspace{-0.51mm}i}+\hspace{-1mm}\Gamma_{\theta_{N\hspace{-0.51mm}-\hspace{-0.51mm}
d\hspace{-0.51mm}-\hspace{-0.51mm}1}}u^{c}_{N\hspace{-0.51mm}-\hspace{-0.51mm}d},\label{f0101}
\end{eqnarray}
and
\begin{eqnarray}
\lambda_{N-1}\hspace{-1mm}&=&\hspace{-1mm}\mathbb{E}_{N-1}\big[Q\hspace{-1mm}+\hspace{-1mm}\bar{A}_{\theta_{N}}'
\bar{P}_{N+1}
\bar{A}_{\theta_{N}}\big]z_{N}\hspace{-1mm}+\hspace{-1mm}\tilde{S}^{1}_{\theta_{N-1}}u^{c}_{N-d}\nonumber\\
\hspace{-1mm}&=&\hspace{-1mm}\bar{P}_{\theta_{N-1}}z_{N}\hspace{-1mm}-\hspace{-1mm}\tilde{S}^{1}_{\theta_{N-1}}\Gamma^{-1}_{\theta_{N-d-1}}
\Bigg(M^{0}_{\theta_{N-d-1}}z_{N-d}\nonumber\\
&&\hspace{-1mm}+\hspace{-1mm}\sum^{d}_{i=1}M^{i}_{\theta_{N-d-1}}u^{c}_{N-d-i}\Bigg),\label{f0102}
\end{eqnarray}
which hold for (\ref{f14}) and (\ref{f16}) in case of $k=N$.\\
Taking $d\leq l\leq N$,  assume that $\Gamma_{\theta_{k-d-1}}(k)>0$, we have (\ref{f14}) and (\ref{f16}) for $k\geq l+1$. Finally, we prove $\Gamma_{\theta_{l-d-1}}(l)>0$.\\
From (\ref{f07}), we have
\begin{eqnarray}
\hspace{-1mm}&&\hspace{-1mm}\mathbb{E}\big[z_{k}'\lambda_{k-1}-z_{k+1}'\lambda_{k}\big]\hspace{-1mm}
=\hspace{-1mm}
\mathbb{E}\big[z_{k}'Qz_{k}
+(u^{c}_{k-d})'Ru^{c}_{k-d}\big].\label{f0103}
\end{eqnarray}
Adding from $k=l+1$ to $k=N$ on both sides of (\ref{f0103}), and taking $z_{l}=0$, we have
\begin{eqnarray}
\tilde{J}_{l}
\hspace{-1mm}&=&\hspace{-1mm}\mathbb{E}\Bigg\{(u^{c}_{l-d})'Ru^{c}_{l-d}
\hspace{-1mm}+\hspace{-1mm}(u^{c}_{l-d})'\bar{B}_{\theta_{l}}'\Bigg[\bar{P}_{\theta_{l}}
\bar{B}_{\theta_{l}}u^{c}_{l-d}\nonumber\\
\hspace{-1mm}&&-\hspace{-1mm}\sum^{d}_{s=0}\left(F^{s+1}_{\theta_{l}}
\Gamma^{-1}_{\theta_{l-d+s}}M^{0}_{\theta_{l-d+s}}
\prod^{s-1}_{i=0}\bar{A}_{\theta_{l-d+i+1}}\right)z_{l-d+1}\nonumber\\
\hspace{-1mm}&&-\hspace{-1mm}\sum^{d}_{i=1}\left(\sum^{d}_{s=0}F^{s+1}_{\theta_{l}}
\Gamma^{-1}_{\theta_{l-d+s}}M^{i+s}_{\theta_{l-d+s}}\right)u^{c}_{l-d-i+1}\Bigg]\Bigg\}\nonumber\\
\hspace{-1mm}&=&\hspace{-1mm}\mathbb{E}\Bigg\{(u^{c}_{l-d})'\mathbb{E}_{l-d-1}\Bigg[R
\hspace{-1mm}+\hspace{-1mm}\bar{B}_{\theta_{l}}'\bar{P}_{\theta_{l}}\bar{B}_{\theta_{l}}
\hspace{-1mm}-\hspace{-1mm}\sum^{d}_{s=0}(M^{s+1}_{\theta_{l-d+s}})'\nonumber\\
\hspace{-1mm}&&\times\Gamma^{-1}_{\theta_{l-d+s}}
M^{s+1}_{\theta_{l-d+s}}\Bigg]u^{c}_{l-d}\Bigg\}\nonumber\\
\hspace{-1mm}&=&\hspace{-1mm}\mathbb{E}\Big[(u^{c}_{l-d})'\Gamma_{\theta_{l-d-1}}
u^{c}_{l-d}\Big].\label{f104}
\end{eqnarray}\\
Due to the uniqueness of the optimal control,  for any nonzero $u(l-d)$, we have $\Gamma_{\theta_{l-d-1}}>0$.\\
It follows from (\ref{f07}) that we have
\begin{small}\begin{eqnarray}
0\hspace{-1mm}&=&\hspace{-1mm}\mathbb{E}_{l\hspace{-0.51mm}-\hspace{-0.51mm}
d\hspace{-0.51mm}-\hspace{-0.51mm}1}\Bigg\{\bar{B}_{\theta_{l}}'
\Bigg[\bar{P}_{\theta_{l}}
(\bar{A}_{\theta_{l}}z_{l}\hspace{-1mm}+\hspace{-1mm}\bar{B}_{\theta_{l}}
u^{c}_{l\hspace{-0.51mm}-\hspace{-0.51mm}d})
\hspace{-1mm}-\hspace{-1mm}\sum^{d}_{s=0}\Big(F^{s\hspace{-0.51mm}+\hspace{-0.51mm}1}_{\theta_{l}}
\Gamma^{-1}_{\theta_{l\hspace{-0.51mm}-\hspace{-0.51mm}d\hspace{-0.51mm}+\hspace{-0.51mm}s}}\nonumber\\
\hspace{-1mm}&&\times\hspace{-1mm}
M^{0}_{\theta_{l\hspace{-0.51mm}-\hspace{-0.51mm}d\hspace{-0.51mm}+\hspace{-0.51mm}s}}
\prod^{s-1}_{i=0}\bar{A}_{\theta_{l\hspace{-0.51mm}-\hspace{-0.51mm}d\hspace{-0.51mm}
+\hspace{-0.51mm}i\hspace{-0.51mm}+\hspace{-0.51mm}1}}\Big)
(\bar{A}_{\theta_{l-d}}z_{l\hspace{-0.51mm}-\hspace{-0.51mm}d}\hspace{-1mm}
+\hspace{-1mm}\bar{B}_{\theta_{l\hspace{-0.51mm}-\hspace{-0.51mm}d}}
u^{c}_{l\hspace{-0.51mm}-\hspace{-0.51mm}2d})\nonumber\\
&&\times
\hspace{-1mm}-\hspace{-1mm}\sum^{d}_{i=1}\left(\sum^{d}_{s=0}
F^{s\hspace{-0.51mm}+\hspace{-0.51mm}1}_{\theta_{l}}
\Gamma^{-1}_{\theta_{l\hspace{-0.51mm}-\hspace{-0.51mm}d\hspace{-0.51mm}+\hspace{-0.51mm}s}}
M^{i\hspace{-0.51mm}+\hspace{-0.51mm}s}_{\theta_{l\hspace{-0.51mm}-\hspace{-0.51mm}d\hspace{-0.51mm}
+\hspace{-0.51mm}s}}\right)u^{c}_{l\hspace{-0.51mm}-\hspace{-0.51mm}d\hspace{-0.51mm}
-\hspace{-0.51mm}i\hspace{-0.51mm}+\hspace{-0.51mm}1}\Bigg]
\hspace{-1mm}+\hspace{-1mm}Ru^{c}_{l\hspace{-0.51mm}-\hspace{-0.51mm}d}\Bigg\}\nonumber\\
\hspace{-1mm}&=&\hspace{-1mm}\mathbb{E}_{l\hspace{-0.51mm}-\hspace{-0.51mm}d\hspace{-0.51mm}-\hspace{-0.51mm}1}
\Bigg\{\Gamma_{\theta_{l\hspace{-0.51mm}-\hspace{-0.51mm}d\hspace{-0.51mm}-\hspace{-0.51mm}1}}
u^{c}_{l\hspace{-0.51mm}-\hspace{-0.51mm}d}\hspace{-1mm}+\hspace{-1mm}
\Gamma_{\theta_{l\hspace{-0.51mm}-\hspace{-0.51mm}d\hspace{-0.51mm}-\hspace{-0.51mm}1}}
u^{c}_{l\hspace{-0.51mm}-\hspace{-0.51mm}d}
\hspace{-1mm}+\hspace{-1mm}(\tilde{S}^{1}_{\theta_{l}})'z_{l}\nonumber\\
&&
\hspace{-1mm}-\hspace{-1mm}\sum^{d}_{s=0}
\Bigg[(M^{s\hspace{-0.51mm}+\hspace{-0.51mm}1}_{\theta_{l\hspace{-0.51mm}-\hspace{-0.51mm}d\hspace{-0.51mm}
+\hspace{-0.51mm}s}})'
\Gamma^{-1}_{\theta_{l\hspace{-0.51mm}-\hspace{-0.51mm}d\hspace{-0.51mm}+\hspace{-0.51mm}s}}
M^{0}_{\theta_{l\hspace{-0.51mm}-\hspace{-0.51mm}d\hspace{-0.51mm}+\hspace{-0.51mm}s}}
\prod^{s}_{i=0}\bar{A}_{\theta_{l\hspace{-0.51mm}-\hspace{-0.51mm}d\hspace{-0.51mm}
+\hspace{-0.51mm}i}}\Bigg]z_{l\hspace{-0.51mm}-\hspace{-0.51mm}d}\nonumber\\
\hspace{-1mm}&&-\hspace{-1mm}\sum^{d}_{i=1}\left(\sum^{d}_{s=0}(M^{s\hspace{-0.51mm}
+\hspace{-0.51mm}1}_{\theta_{l\hspace{-0.51mm}-\hspace{-0.51mm}d\hspace{-0.51mm}+\hspace{-0.51mm}s}})'
\Gamma^{-1}_{\theta_{l\hspace{-0.51mm}-\hspace{-0.51mm}d\hspace{-0.51mm}+\hspace{-0.51mm}s}}
M^{i\hspace{-0.51mm}+\hspace{-0.51mm}s\hspace{-0.51mm}+\hspace{-0.51mm}1}_{\theta_{l\hspace{-0.51mm}
-\hspace{-0.51mm}d\hspace{-0.51mm}+\hspace{-0.51mm}s}}\right)u^{c}_{l\hspace{-0.51mm}
-\hspace{-0.51mm}d\hspace{-0.51mm}-\hspace{-0.51mm}i}\Bigg\}
\nonumber\\
\hspace{-1mm}&=&\hspace{-1mm}\Gamma_{\theta_{l\hspace{-0.51mm}-\hspace{-0.51mm}d\hspace{-0.51mm}
-\hspace{-0.51mm}1}}u^{c}_{l\hspace{-0.51mm}-\hspace{-0.51mm}d}\hspace{-1mm}+\hspace{-1mm}
\mathbb{E}_{l\hspace{-0.51mm}-\hspace{-0.51mm}d\hspace{-0.51mm}-\hspace{-0.51mm}1}
\Bigg[(\tilde{S}^{1}_{\theta_{l\hspace{-0.51mm}-\hspace{-0.51mm}1}})'
\prod^{d}_{j=1}\bar{A}_{\theta_{l\hspace{-0.51mm}-\hspace{-0.51mm}j}}\hspace{-1mm}-\hspace{-1mm}
\sum^{d}_{i=0}\Bigg((M^{i\hspace{-0.51mm}+\hspace{-0.51mm}1}_{\theta_{l\hspace{-0.51mm}-\hspace{-0.51mm}
d\hspace{-0.51mm}+\hspace{-0.51mm}i}})'\nonumber\\
\hspace{-1mm}&&\times\hspace{-1mm}\Gamma^{-1}_{\theta_{l\hspace{-0.51mm}-\hspace{-0.51mm}d\hspace{-0.51mm}+\hspace{-0.51mm}i}}
M^{0}_{\theta_{l\hspace{-0.51mm}-\hspace{-0.51mm}d\hspace{-0.51mm}+\hspace{-0.51mm}i}}
\prod^{i}_{s=0}
\bar{A}_{\theta_{l\hspace{-0.51mm}-\hspace{-0.51mm}d\hspace{-0.51mm}+\hspace{-0.51mm}s}}\Bigg)\Bigg]
z_{l\hspace{-0.51mm}-\hspace{-0.51mm}d}
\hspace{-1mm}-\hspace{-1mm}\mathbb{E}_{l\hspace{-0.51mm}-\hspace{-0.51mm}d\hspace{-0.51mm}
-\hspace{-0.51mm}1}\Bigg[(\tilde{S}^{1}_{\theta_{l\hspace{-0.51mm}-\hspace{-0.51mm}1}})'
\nonumber\\
\hspace{-1mm}&&\times\hspace{-1mm}\prod^{i-1}_{j=1}\bar{A}_{\theta_{l\hspace{-0.51mm}-\hspace{-0.51mm}j}}
\bar{B}_{\theta_{l\hspace{-0.51mm}
-\hspace{-0.51mm}i}}\hspace{-1mm}-\hspace{-1mm}\sum^{d}_{s=0}
\Big((M^{s\hspace{-0.51mm}+\hspace{-0.51mm}1}_{\theta_{l\hspace{-0.51mm}-\hspace{-0.51mm}d\hspace{-0.51mm}
+\hspace{-0.51mm}s}})'
\Gamma^{-1}_{\theta_{l\hspace{-0.51mm}-\hspace{-0.51mm}d\hspace{-0.51mm}+\hspace{-0.51mm}s}}
M^{i\hspace{-0.51mm}+\hspace{-0.51mm}s\hspace{-0.51mm}+\hspace{-0.51mm}1}_{\theta_{l\hspace{-0.51mm}
-\hspace{-0.51mm}d\hspace{-0.51mm}+\hspace{-0.51mm}s}}\Big)\Bigg]u^{c}_{l\hspace{-0.51mm}-\hspace{-0.51mm}
d\hspace{-0.51mm}-\hspace{-0.51mm}i},\label{f105}
\end{eqnarray}
\end{small}
and
\begin{small}
\begin{eqnarray}
\lambda_{l-1}\hspace{-1mm}&=&\hspace{-1mm}\mathbb{E}_{l\hspace{-0.51mm}-\hspace{-0.51mm}1}
\Bigg\{Qz_{l}
\hspace{-1mm}+\hspace{-1mm}\bar{A}_{\theta_{l}}'\Big[\bar{P}_{\theta_{l}}
(\bar{A}_{\theta_{l}}z_{l}\hspace{-1mm}+\hspace{-1mm}\bar{B}_{\theta_{l}}
u^{c}_{l\hspace{-0.51mm}-\hspace{-0.51mm}d})
\hspace{-1mm}-\hspace{-1mm}\sum^{d}_{s=0}\Bigg(F^{s\hspace{-0.51mm}
+\hspace{-0.51mm}1}_{\theta_{l}}
\Gamma^{-1}_{\theta_{l\hspace{-0.51mm}-\hspace{-0.51mm}d\hspace{-0.51mm}+\hspace{-0.51mm}s}}\nonumber\\
\hspace{-1mm}&&\times\hspace{-1mm}
M^{0}_{\theta_{l\hspace{-0.51mm}-\hspace{-0.51mm}d\hspace{-0.51mm}+\hspace{-0.51mm}s}}
\prod^{s-1}_{i=0}\bar{A}_{\theta_{l\hspace{-0.51mm}-\hspace{-0.51mm}d\hspace{-0.51mm}
+\hspace{-0.51mm}i\hspace{-0.51mm}+\hspace{-0.51mm}1}}\Bigg)
(\bar{A}_{\theta_{l-d}}z_{l-d}\hspace{-1mm}+\hspace{-1mm}\bar{B}_{\theta_{l\hspace{-0.51mm}
-\hspace{-0.51mm}d}}u^{c}_{l\hspace{-0.51mm}-\hspace{-0.51mm}2d})\nonumber\\
\hspace{-1mm}&&-\hspace{-1mm}\sum^{d}_{i=1}\left(\sum^{d}_{s=0}F^{s\hspace{-0.51mm}+\hspace{-0.51mm}1}_{\theta_{l}}
\Gamma^{-1}_{\theta_{l\hspace{-0.51mm}-\hspace{-0.51mm}d\hspace{-0.51mm}+\hspace{-0.51mm}s}}
M^{i\hspace{-0.51mm}+\hspace{-0.51mm}s}_{\theta_{l\hspace{-0.51mm}-\hspace{-0.51mm}
d\hspace{-0.51mm}+\hspace{-0.51mm}s}}\right)u^{c}_{l\hspace{-0.51mm}-\hspace{-0.51mm}
d\hspace{-0.51mm}-\hspace{-0.51mm}i\hspace{-0.51mm}+\hspace{-0.51mm}1}\Big]\Bigg\}\nonumber\\
\hspace{-1mm}&=&\hspace{-1mm}\bar{P}_{\theta_{l\hspace{-0.51mm}-\hspace{-0.51mm}1}}z_{l}
\hspace{-1mm}+\hspace{-1mm}\Bigg\{\tilde{S}^{1}_{\theta_{l\hspace{-0.51mm}-\hspace{-0.51mm}1}}\hspace{-1mm}-\hspace{-1mm}
\sum^{d}_{s=0}\Big[\mathbb{E}_{l\hspace{-0.51mm}-\hspace{-0.51mm}1}\bar{A}_{\theta_{l}}'F^{s\hspace{-0.51mm}+\hspace{-0.51mm}1}_{\theta_{l}}
\Gamma^{-1}_{\theta_{l\hspace{-0.51mm}-\hspace{-0.51mm}d\hspace{-0.51mm}+\hspace{-0.51mm}s}}M^{s\hspace{-0.51mm}+\hspace{-0.51mm}1}_{\theta_{l\hspace{-0.51mm}-\hspace{-0.51mm}d\hspace{-0.51mm}+\hspace{-0.51mm}s}}\Big]\Bigg\}u^{c}_{l\hspace{-0.51mm}-\hspace{-0.51mm}d}\nonumber\\
\hspace{-1mm}&&-\hspace{-1mm}\mathbb{E}_{l\hspace{-0.51mm}-\hspace{-0.51mm}1}\Bigg\{\sum^{d}_{s=0}
\Bigg(\bar{A}_{\theta_{l}}'F^{s\hspace{-0.51mm}+\hspace{-0.51mm}1}_{\theta_{l}}
\Gamma^{-1}_{\theta_{l\hspace{-0.51mm}-\hspace{-0.51mm}d\hspace{-0.51mm}+\hspace{-0.51mm}s}}M^{0}_{\theta_{l\hspace{-0.51mm}-\hspace{-0.51mm}d\hspace{-0.51mm}+\hspace{-0.51mm}s}}\prod^{s}_{i=0}\bar{A}_{\theta_{l\hspace{-0.51mm}-\hspace{-0.51mm}d\hspace{-0.51mm}+\hspace{-0.51mm}i}}\Bigg)
z_{l\hspace{-0.51mm}-\hspace{-0.51mm}d}\nonumber\\
\hspace{-1mm}&&
+\hspace{-1mm}\sum^{d}_{i=1}\left(\sum^{d}_{s=0}\bar{A}_{\theta_{l}}'F^{s\hspace{-0.51mm}+\hspace{-0.51mm}1}_{\theta_{l}}
\Gamma^{-1}_{\theta_{l\hspace{-0.51mm}-\hspace{-0.51mm}d\hspace{-0.51mm}+\hspace{-0.51mm}s}}
M^{i\hspace{-0.51mm}+\hspace{-0.51mm}s\hspace{-0.51mm}+\hspace{-0.51mm}1}_{\theta_{l
\hspace{-0.51mm}-\hspace{-0.51mm}d\hspace{-0.51mm}+\hspace{-0.51mm}s}}\right)
u^{c}_{l\hspace{-0.51mm}-\hspace{-0.51mm}d\hspace{-0.51mm}-\hspace{-0.51mm}i}\Big]\Bigg\}\nonumber\\
\hspace{-1mm}&=&\hspace{-1mm}\bar{P}_{\theta_{l\hspace{-0.51mm}-\hspace{-0.51mm}1}}z_{l}
\hspace{-1mm}+\hspace{-1mm}\Bigg[\tilde{S}^{1}_{\theta_{l\hspace{-0.51mm}-\hspace{-0.51mm}1}}\hspace{-1mm}-\hspace{-1mm}
\sum^{d}_{s=0}\Big(F^{s\hspace{-0.51mm}+\hspace{-0.51mm}2}_{\theta_{l\hspace{-0.51mm}-\hspace{-0.51mm}1}}
\Gamma^{-1}_{\theta_{l\hspace{-0.51mm}-\hspace{-0.51mm}d\hspace{-0.51mm}+\hspace{-0.51mm}s}}
M^{s\hspace{-0.51mm}+\hspace{-0.51mm}1}_{\theta_{l\hspace{-0.51mm}-\hspace{-0.51mm}d
\hspace{-0.51mm}+\hspace{-0.51mm}s}}\Big)\Bigg]\nonumber\\
\hspace{-1mm}&&-\hspace{-1mm}\Bigg\{\sum^{d}_{s=0}\Big(F^{s\hspace{-0.51mm}+\hspace{-0.51mm}2}_{\theta_{l\hspace{-0.51mm}-\hspace{-0.51mm}1}}
\Gamma^{-1}_{\theta_{l\hspace{-0.51mm}-\hspace{-0.51mm}d\hspace{-0.51mm}+\hspace{-0.51mm}s}}M^{0}_{\theta_{l\hspace{-0.51mm}-\hspace{-0.51mm}d\hspace{-0.51mm}+\hspace{-0.51mm}s}}\prod^{s}_{i=0}\bar{A}_{\theta_{l\hspace{-0.51mm}-\hspace{-0.51mm}d\hspace{-0.51mm}+\hspace{-0.51mm}i}}\Big)
z_{l\hspace{-0.51mm}-\hspace{-0.51mm}d}\nonumber
\end{eqnarray}
\end{small}
\begin{small}
\begin{eqnarray}
\hspace{-1mm}&&+\hspace{-1mm}\sum^{d}_{i=1}\left(\sum^{d}_{s=0}F^{s\hspace{-0.51mm}+\hspace{-0.51mm}2}_{\theta_{l\hspace{-0.51mm}-\hspace{-0.51mm}1}}
\Gamma^{-1}_{\theta_{l\hspace{-0.51mm}-\hspace{-0.51mm}d\hspace{-0.51mm}+\hspace{-0.51mm}s}}
M^{i\hspace{-0.51mm}+\hspace{-0.51mm}s\hspace{-0.51mm}+\hspace{-0.51mm}1}_{\theta_{l
\hspace{-0.51mm}-\hspace{-0.51mm}d\hspace{-0.51mm}+\hspace{-0.51mm}s}}\right)
u^{c}_{l\hspace{-0.51mm}-\hspace{-0.51mm}d\hspace{-0.51mm}-\hspace{-0.51mm}i}\Big]\Bigg\},\label{f106}
\end{eqnarray}
\end{small}
i.e., $u^{c}_{l-d}$ and $\lambda_{l-1}$ are as (\ref{f14}) and  (\ref{f16}) with $k=l$, respectively.\\
($\Longleftarrow$) When $\Gamma_{\theta_{k-d-1}}(k)>0$, we will investigate the unique solvability of Problem 2. Define
\begin{small}
\begin{eqnarray}
V_{N}(k)\hspace{-1mm}&=&\hspace{-1mm}\mathbb{E}\Bigg\{z_{k}'\bar{P}_{\theta_{k\hspace{-0.51mm}
-\hspace{-0.51mm}1}}z_{k}
\hspace{-1mm}-\hspace{-1mm}z_{k}'\sum^{d}_{s=0}\left(F^{s\hspace{-0.51mm}+\hspace{-0.51mm}1}_{\theta_{k\hspace{-0.51mm}-\hspace{-0.51mm}1}}
\Gamma^{-1}_{\theta_{k\hspace{-0.51mm}-\hspace{-0.51mm}d\hspace{-0.51mm}-\hspace{-0.51mm}1\hspace{-0.51mm}+\hspace{-0.51mm}s}}M^{0}_{\theta_{k\hspace{-0.51mm}-\hspace{-0.51mm}d\hspace{-0.51mm}-\hspace{-0.51mm}1\hspace{-0.51mm}+\hspace{-0.51mm}s}}
\prod^{s-1}_{i=0}\bar{A}_{\theta_{k\hspace{-0.51mm}-\hspace{-0.51mm}d\hspace{-0.51mm}
+\hspace{-0.51mm}i}}\right)\nonumber\\
\hspace{-1mm}&&\times\hspace{-1mm}z_{k\hspace{-0.51mm}-\hspace{-0.51mm}d}\hspace{-1mm}-\hspace{-1mm}z_{k}'\sum^{d}_{i=1}\left(\sum^{d}_{s=0}F^{s\hspace{-0.51mm}+\hspace{-0.51mm}1}_{\theta_{k\hspace{-0.51mm}-\hspace{-0.51mm}1}}
\Gamma^{-1}_{\theta_{k\hspace{-0.51mm}-\hspace{-0.51mm}d\hspace{-0.51mm}-\hspace{-0.51mm}1\hspace{-0.51mm}+\hspace{-0.51mm}s}}M^{i\hspace{-0.51mm}+\hspace{-0.51mm}s}_{\theta_{k\hspace{-0.51mm}-\hspace{-0.51mm}d\hspace{-0.51mm}-\hspace{-0.51mm}1\hspace{-0.51mm}+\hspace{-0.51mm}s}}\right)u^{c}_{k\hspace{-0.51mm}-\hspace{-0.51mm}d\hspace{-0.51mm}-\hspace{-0.51mm}i}\Bigg\}. \label{f107}
\end{eqnarray}
\end{small}
From Remark 4, we have
\begin{small}
\begin{eqnarray}
\hspace{-1mm}&&\hspace{-1mm}V_{N}(k)-V_{N}(k+1)\nonumber\\
\hspace{-1mm}&=&\hspace{-1mm}
\mathbb{E}\Bigg\{z_{k}'[\bar{P}_{\theta_{k\hspace{-0.51mm}-\hspace{-0.51mm}1}}\hspace{-1mm}-\hspace{-1mm}
\bar{A}_{\theta_{k}}'\bar{P}_{\theta_{k}}
\bar{A}_{\theta_{k}}]z_{k}
\hspace{-1mm}-\hspace{-1mm}z_{k}'\sum^{d}_{s=0}\Bigg[\bar{A}_{\theta_{k}}'\bar{P}_{\theta_{k}}
\bar{B}_{\theta_{k}}
\hspace{-1mm}-\hspace{-1mm}F^{s\hspace{-0.51mm}+\hspace{-0.51mm}2}_{\theta_{k\hspace{-0.51mm}-\hspace{-0.51mm}1}}
\nonumber\\
\hspace{-1mm}&&\times\hspace{-1mm}
\Gamma^{-1}_{\theta_{k\hspace{-0.51mm}-\hspace{-0.51mm}d\hspace{-0.51mm}+\hspace{-0.51mm}s}}
M^{s\hspace{-0.51mm}+\hspace{-0.51mm}1}_{\theta_{k\hspace{-0.51mm}-\hspace{-0.51mm}d
\hspace{-0.51mm}+\hspace{-0.51mm}s}}\Bigg]u^{c}_{k\hspace{-0.51mm}-\hspace{-0.51mm}d}
-\hspace{-1mm}(u^{c}_{k\hspace{-0.51mm}-\hspace{-0.51mm}d})'\bar{B}_{\theta_{k}}'\bar{P}_{\theta_{k}}
\bar{A}_{\theta_{k}}z_{k}\nonumber\\
\hspace{-1mm}&&-\hspace{-1mm}z_{k}'\Bigg[\sum^{d}_{s=0}(F^{s\hspace{-0.51mm}+\hspace{-0.51mm}1}_{\theta_{k\hspace{-0.51mm}-\hspace{-0.51mm}1}}
\Gamma^{-1}_{\theta_{k\hspace{-0.51mm}-\hspace{-0.51mm}d\hspace{-0.51mm}-\hspace{-0.51mm}1\hspace{-0.51mm}+\hspace{-0.51mm}s}}M^{0}_{\theta_{k\hspace{-0.51mm}-\hspace{-0.51mm}d\hspace{-0.51mm}-\hspace{-0.51mm}1\hspace{-0.51mm}+\hspace{-0.51mm}s}}
\prod^{s-1}_{i=0}\bar{A}_{\theta_{k\hspace{-0.51mm}-\hspace{-0.51mm}d\hspace{-0.51mm}+\hspace{-0.51mm}i}})\nonumber\\
\hspace{-1mm}&&-\hspace{-1mm}\sum^{d}_{s=0}(F^{s\hspace{-0.51mm}+\hspace{-0.51mm}2}_{\theta_{k\hspace{-0.51mm}-\hspace{-0.51mm}1}}
\Gamma^{-1}_{\theta_{k\hspace{-0.51mm}-\hspace{-0.51mm}d\hspace{-0.51mm}+\hspace{-0.51mm}s}}M^{0}_{\theta_{k\hspace{-0.51mm}-\hspace{-0.51mm}d\hspace{-0.51mm}+\hspace{-0.51mm}s}}
\prod^{s}_{i=0}\bar{A}_{\theta_{k\hspace{-0.51mm}-\hspace{-0.51mm}d\hspace{-0.51mm}
+\hspace{-0.51mm}i}})\Bigg]z_{k\hspace{-0.51mm}-\hspace{-0.51mm}d}
\hspace{-1mm}-\hspace{-1mm}(u^{c}_{k\hspace{-0.51mm}-\hspace{-0.51mm}d})'\nonumber\\
\hspace{-1mm}&&\times\hspace{-1mm}\Bigg[\bar{B}_{\theta_{k}}'\bar{P}_{\theta_{k}}
\bar{B}_{\theta_{k}}
\hspace{-1mm}-\hspace{-1mm}\sum^{d}_{s=0}(M^{s\hspace{-0.51mm}+\hspace{-0.51mm}1}_{\theta_{k\hspace{-0.51mm}-\hspace{-0.51mm}d\hspace{-0.51mm}+\hspace{-0.51mm}s}})'
\Gamma^{-1}_{\theta_{k\hspace{-0.51mm}-\hspace{-0.51mm}d\hspace{-0.51mm}+\hspace{-0.51mm}s}}M^{s\hspace{-0.51mm}+\hspace{-0.51mm}1}_{\theta_{k\hspace{-0.51mm}-\hspace{-0.51mm}d\hspace{-0.51mm}+\hspace{-0.51mm}s}}\Bigg]u^{c}_{k\hspace{-0.51mm}-\hspace{-0.51mm}d}\nonumber\\
\hspace{-1mm}&&+\hspace{-1mm}(u^{c}_{k\hspace{-0.51mm}-\hspace{-0.51mm}d})'\sum^{d}_{s=0}\left((M^{s\hspace{-0.51mm}+\hspace{-0.51mm}1}_{\theta_{k\hspace{-0.51mm}-\hspace{-0.51mm}d\hspace{-0.51mm}+\hspace{-0.51mm}s}})'
\Gamma^{-1}_{\theta_{k\hspace{-0.51mm}-\hspace{-0.51mm}d\hspace{-0.51mm}+\hspace{-0.51mm}s}}M^{0}_{\theta_{k\hspace{-0.51mm}-\hspace{-0.51mm}d\hspace{-0.51mm}+\hspace{-0.51mm}s}}
\prod^{s}_{i=0}\bar{A}_{\theta_{k\hspace{-0.51mm}-\hspace{-0.51mm}d\hspace{-0.51mm}+\hspace{-0.51mm}i}}\right)z_{k\hspace{-0.51mm}-\hspace{-0.51mm}d}\nonumber\\
\hspace{-1mm}&&-\hspace{-1mm}z_{k}'\sum^{d}_{i=1}\left(F^{s\hspace{-0.51mm}+\hspace{-0.51mm}1}_{\theta_{k\hspace{-0.51mm}-\hspace{-0.51mm}1}}
\Gamma^{-1}_{\theta_{k\hspace{-0.51mm}-\hspace{-0.51mm}d\hspace{-0.51mm}+\hspace{-0.51mm}s\hspace{-0.51mm}-\hspace{-0.51mm}1}}M^{i\hspace{-0.51mm}+\hspace{-0.51mm}s}_{\theta_{k\hspace{-0.51mm}-\hspace{-0.51mm}d\hspace{-0.51mm}+\hspace{-0.51mm}s\hspace{-0.51mm}-\hspace{-0.51mm}1}}-F^{s\hspace{-0.51mm}+\hspace{-0.51mm}2}_{\theta_{k\hspace{-0.51mm}-\hspace{-0.51mm}1}}
\Gamma^{-1}_{\theta_{k\hspace{-0.51mm}-\hspace{-0.51mm}d\hspace{-0.51mm}+\hspace{-0.51mm}s}}M^{i\hspace{-0.51mm}+\hspace{-0.51mm}s\hspace{-0.51mm}+\hspace{-0.51mm}1}_{\theta_{k\hspace{-0.51mm}-\hspace{-0.51mm}d\hspace{-0.51mm}+\hspace{-0.51mm}s}}\right)u^{c}_{k\hspace{-0.51mm}-\hspace{-0.51mm}d\hspace{-0.51mm}-\hspace{-0.51mm}i}\nonumber\\
\hspace{-1mm}&&+\hspace{-1mm}(u^{c}_{k\hspace{-0.51mm}-\hspace{-0.51mm}d})'\sum^{d}_{i=1}\left(\sum^{d}_{s=0}
(M^{s\hspace{-0.51mm}+\hspace{-0.51mm}1}_{\theta_{k\hspace{-0.51mm}-\hspace{-0.51mm}d\hspace{-0.51mm}+\hspace{-0.51mm}s}})'
\Gamma^{-1}_{\theta_{k\hspace{-0.51mm}-\hspace{-0.51mm}d\hspace{-0.51mm}+\hspace{-0.51mm}s}}
M^{i\hspace{-0.51mm}+\hspace{-0.51mm}s\hspace{-0.51mm}+\hspace{-0.51mm}1}_{\theta_{k
\hspace{-0.51mm}-\hspace{-0.51mm}d\hspace{-0.51mm}+\hspace{-0.51mm}s}}\right)
u^{c}_{k\hspace{-0.51mm}-\hspace{-0.51mm}d\hspace{-0.51mm}-\hspace{-0.51mm}i}\Bigg\}\nonumber\\
\hspace{-1mm}&=&\hspace{-1mm}
\mathbb{E}\Bigg\{z_{k}'Qz_{k}-(u^{c}_{k\hspace{-0.51mm}-\hspace{-0.51mm}d})'\big(
\Gamma_{\theta_{k\hspace{-0.51mm}-\hspace{-0.51mm}d\hspace{-0.51mm}-\hspace{-0.51mm}1}}
-R\big)u^{c}_{k\hspace{-0.51mm}-\hspace{-0.51mm}d}
\hspace{-1mm}-\hspace{-1mm}z_{k}'F^{1}_{\theta_{k\hspace{-0.51mm}-\hspace{-0.51mm}1}}
u^{c}_{k\hspace{-0.51mm}-\hspace{-0.51mm}d}\nonumber\\
\hspace{-1mm}&&-\hspace{-1mm}
z_{k}'F^{1}_{\theta_{k\hspace{-0.51mm}-\hspace{-0.51mm}1}}
\hspace{-1mm}-\hspace{-1mm}z_{k}'\Bigg(F^{1}_{\theta_{k\hspace{-0.51mm}-\hspace{-0.51mm}1}}
\Gamma^{-1}_{\theta_{k\hspace{-0.51mm}-\hspace{-0.51mm}d\hspace{-0.51mm}-\hspace{-0.51mm}1}}
M^{0}_{\theta_{k\hspace{-0.51mm}-\hspace{-0.51mm}d\hspace{-0.51mm}-\hspace{-0.51mm}1}}
\hspace{-1mm}-\hspace{-1mm}F^{d\hspace{-0.51mm}+\hspace{-0.51mm}2}_{\theta_{k\hspace{-0.51mm}
-\hspace{-0.51mm}1}}\Gamma^{-1}_{\theta_{k}}\hspace{-1mm}+\hspace{-1mm}
(u^{c}_{k\hspace{-0.51mm}-\hspace{-0.51mm}d})'
\nonumber\\
\hspace{-1mm}&&\times\hspace{-1mm}\sum^{d}_{s=0}\Bigg((M^{s\hspace{-0.51mm}
+\hspace{-0.51mm}1}_{\theta_{k\hspace{-0.51mm}-\hspace{-0.51mm}d\hspace{-0.51mm}+\hspace{-0.51mm}s}})'
\Gamma^{-1}_{\theta_{k\hspace{-0.51mm}-\hspace{-0.51mm}d\hspace{-0.51mm}+\hspace{-0.51mm}s}}
M^{0}_{\theta_{k}}\prod^{d}_{j=0}\bar{A}_{\theta_{k\hspace{-0.51mm}-\hspace{-0.51mm}d
\hspace{-0.51mm}+\hspace{-0.51mm}j}}\Bigg)z_{k\hspace{-0.51mm}-\hspace{-0.51mm}d}\nonumber\\
\hspace{-1mm}&&\times \hspace{-1mm}
M^{0}_{\theta_{k\hspace{-0.51mm}-\hspace{-0.51mm}d\hspace{-0.51mm}+\hspace{-0.51mm}s}}
\prod^{s}_{i=0}\bar{A}_{\theta_{k\hspace{-0.51mm}-\hspace{-0.51mm}d\hspace{-0.51mm}
+\hspace{-0.51mm}i}}\Bigg)z_{k\hspace{-0.51mm}-\hspace{-0.51mm}d}\hspace{-1mm}-\hspace{-1mm}
z_{k}'\sum^{d}_{i=1}\sum^{d}_{s=0}\Bigg(F^{1}_{\theta_{k\hspace{-0.51mm}-\hspace{-0.51mm}1}}
\nonumber
\end{eqnarray}
\end{small}
\begin{small}
\begin{eqnarray}
\hspace{-1mm}&&\times \hspace{-1mm}
\Gamma^{-1}_{\theta_{k\hspace{-0.51mm}-\hspace{-0.51mm}d\hspace{-0.51mm}+\hspace{-0.51mm}s\hspace{-0.51mm}-\hspace{-0.51mm}1}}M^{i}_{\theta_{k\hspace{-0.51mm}-\hspace{-0.51mm}d\hspace{-0.51mm}-\hspace{-0.51mm}1}}\hspace{-1mm}-\hspace{-1mm}
F^{d\hspace{-0.51mm}+\hspace{-0.51mm}2}_{\theta_{k\hspace{-0.51mm}-\hspace{-0.51mm}1}}
\Gamma^{-1}_{\theta_{k}}M^{i\hspace{-0.51mm}+\hspace{-0.51mm}d\hspace{-0.51mm}
+\hspace{-0.51mm}1}_{\theta_{k}}\Bigg)u^{c}_{k\hspace{-0.51mm}-\hspace{-0.51mm}d
\hspace{-0.51mm}-\hspace{-0.51mm}i}\nonumber\\
\hspace{-1mm}&&+\hspace{-1mm}(u^{c}_{k\hspace{-0.51mm}-\hspace{-0.51mm}d})'\sum^{d}_{i=1}
\Bigg(\sum^{d}_{s=0}(M^{s\hspace{-0.51mm}+\hspace{-0.51mm}1}_{\theta_{k\hspace{-0.51mm}
-\hspace{-0.51mm}d\hspace{-0.51mm}+\hspace{-0.51mm}s}})'\Gamma^{-1}_{\theta_{k\hspace{-0.51mm}
-\hspace{-0.51mm}d\hspace{-0.51mm}+\hspace{-0.51mm}s}}M^{i\hspace{-0.51mm}+\hspace{-0.51mm}s
\hspace{-0.51mm}+\hspace{-0.51mm}1}_{\theta_{k\hspace{-0.51mm}-\hspace{-0.51mm}d\hspace{-0.51mm}
+\hspace{-0.51mm}s}}\Bigg)\nonumber\\
\hspace{-1mm}&=&\hspace{-1mm}
\mathbb{E}\Bigg\{z_{k}'Qz_{k}\hspace{-1mm}+\hspace{-1mm}(u^{c}_{k\hspace{-0.51mm}-\hspace{-0.51mm}d})'Ru^{c}_{k\hspace{-0.51mm}-\hspace{-0.51mm}d}
\hspace{-1mm}-\hspace{-1mm}(u^{c}_{k\hspace{-0.51mm}-\hspace{-0.51mm}d})'\Gamma_{\theta_{k\hspace{-0.51mm}-\hspace{-0.51mm}d\hspace{-0.51mm}-\hspace{-0.51mm}1}}u^{c}_{k\hspace{-0.51mm}-\hspace{-0.51mm}d}\nonumber\\
\hspace{-1mm}&&-\hspace{-1mm}\Bigg[\prod^{d}_{j=0}\bar{A}_{\theta_{k\hspace{-0.51mm}-\hspace{-0.51mm}j}}z_{k\hspace{-0.51mm}-\hspace{-0.51mm}d}
\hspace{-1mm}+\hspace{-1mm}\sum^{d}_{i=1}\Bigg(\prod^{i-1}_{j=0}\bar{A}_{\theta_{k\hspace{-0.51mm}-\hspace{-0.51mm}j}}
\bar{B}_{\theta_{k\hspace{-0.51mm}-\hspace{-0.51mm}i}}u^{c}_{k\hspace{-0.51mm}-\hspace{-0.51mm}d
\hspace{-0.51mm}-\hspace{-0.51mm}i}\Bigg)\Bigg]'\nonumber\\ \hspace{-1mm}&&\times\hspace{-1mm}
F^{1}_{\theta_{k\hspace{-0.51mm}-\hspace{-0.51mm}1}}u^{c}_{k\hspace{-0.51mm}-\hspace{-0.51mm}d}
\hspace{-1mm}-\hspace{-1mm}\Bigg[\prod^{d}_{j=0}\bar{A}_{\theta_{k\hspace{-0.51mm}-\hspace{-0.51mm}j}}z_{k\hspace{-0.51mm}-\hspace{-0.51mm}d}
\hspace{-1mm}+\hspace{-1mm}\sum^{d}_{i=1}\Bigg(\prod^{i-1}_{j=0}\bar{A}_{\theta_{k\hspace{-0.51mm}-\hspace{-0.51mm}j}}
\bar{B}_{\theta_{k\hspace{-0.51mm}-\hspace{-0.51mm}i}}\nonumber\\ \hspace{-1mm}&&\times\hspace{-1mm}
u^{c}_{k\hspace{-0.51mm}-\hspace{-0.51mm}d\hspace{-0.51mm}-\hspace{-0.51mm}i}\Bigg)\Bigg]'F^{1}_{\theta_{k\hspace{-0.51mm}-\hspace{-0.51mm}1}}
\Gamma^{-1}_{\theta_{k\hspace{-0.51mm}-\hspace{-0.51mm}d\hspace{-0.51mm}-\hspace{-0.51mm}1}}
M^{0}_{\theta_{k\hspace{-0.51mm}-\hspace{-0.51mm}d\hspace{-0.51mm}-\hspace{-0.51mm}1}}
z_{k\hspace{-0.51mm}-\hspace{-0.51mm}d}\hspace{-1mm}+\hspace{-1mm}(u^{c}_{k\hspace{-0.51mm}
-\hspace{-0.51mm}d})'\sum^{d}_{s=0}\Bigg((M^{s\hspace{-0.51mm}+\hspace{-0.51mm}1}_{\theta_{k
\hspace{-0.51mm}-\hspace{-0.51mm}d\hspace{-0.51mm}+\hspace{-0.51mm}s}})'\nonumber\\
\hspace{-1mm}&&+\hspace{-1mm}(u^{c}_{k\hspace{-0.51mm}-\hspace{-0.51mm}d})'\sum^{d}_{s=0}\Bigg((M^{s\hspace{-0.51mm}+\hspace{-0.51mm}1}_{\theta_{k\hspace{-0.51mm}-\hspace{-0.51mm}d\hspace{-0.51mm}+\hspace{-0.51mm}s}})'
\Gamma^{-1}_{\theta_{k\hspace{-0.51mm}-\hspace{-0.51mm}d\hspace{-0.51mm}+\hspace{-0.51mm}s}}
M^{0}_{\theta_{k\hspace{-0.51mm}-\hspace{-0.51mm}d\hspace{-0.51mm}+\hspace{-0.51mm}s}}
\prod^{s}_{j=0}\bar{A}_{\theta_{k\hspace{-0.51mm}-\hspace{-0.51mm}d\hspace{-0.51mm}
+\hspace{-0.51mm}j}}\Bigg)\nonumber\\
\hspace{-1mm}&&\times\hspace{-1mm}
z_{k-d}\hspace{-1mm}-\hspace{-1mm}\Bigg(\prod^{d}_{j=0}\bar{A}_{\theta_{k\hspace{-0.51mm}
-\hspace{-0.51mm}j}}z_{k\hspace{-0.51mm}-\hspace{-0.51mm}d}\hspace{-1mm}+\hspace{-1mm}\sum^{d}_{i=1}\prod^{i-1}_{j=0}\bar{A}_{\theta_{k\hspace{-0.51mm}-\hspace{-0.51mm}j}}
\bar{B}_{\theta_{k\hspace{-0.51mm}-\hspace{-0.51mm}i}}u^{c}_{k\hspace{-0.51mm}-\hspace{-0.51mm}d
\hspace{-0.51mm}-\hspace{-0.51mm}i}\Bigg)'\nonumber\\
\hspace{-1mm}&&\times\hspace{-1mm}
\sum^{d}_{i=1}\Big(F^{1}_{\theta_{k\hspace{-0.51mm}-\hspace{-0.51mm}1}}
\Gamma^{-1}_{\theta_{k\hspace{-0.51mm}-\hspace{-0.51mm}d\hspace{-0.51mm}-\hspace{-0.51mm}1}}
M^{i}_{\theta_{k\hspace{-0.51mm}-\hspace{-0.51mm}d\hspace{-0.51mm}-\hspace{-0.51mm}1}}
u^{c}_{k\hspace{-0.51mm}-\hspace{-0.51mm}d\hspace{-0.51mm}-\hspace{-0.51mm}i}\Big)
\hspace{-1mm}+\hspace{-1mm}(u^{c}_{k\hspace{-0.51mm}-\hspace{-0.51mm}d})'\nonumber\\
\hspace{-1mm}&&\times\hspace{-1mm}\sum^{d}_{i=1}\left(\sum^{d}_{s=0}
(M^{s\hspace{-0.51mm}+\hspace{-0.51mm}1}_{\theta_{k\hspace{-0.51mm}-\hspace{-0.51mm}d\hspace{-0.51mm}+\hspace{-0.51mm}s}})'
\Gamma^{-1}_{\theta_{k\hspace{-0.51mm}-\hspace{-0.51mm}d\hspace{-0.51mm}+\hspace{-0.51mm}s}}
M^{i\hspace{-0.51mm}+\hspace{-0.51mm}s\hspace{-0.51mm}+\hspace{-0.51mm}1}_{\theta_{k
\hspace{-0.51mm}-\hspace{-0.51mm}d\hspace{-0.51mm}+\hspace{-0.51mm}s}}\right)
u^{c}_{k\hspace{-0.51mm}-\hspace{-0.51mm}d\hspace{-0.51mm}-\hspace{-0.51mm}i}\hspace{-1mm}
-\hspace{-1mm}(u^{c}_{k\hspace{-0.51mm}-\hspace{-0.51mm}d})'
(\tilde{S}^{1}_{_{\theta_{k\hspace{-0.51mm}-\hspace{-0.51mm}1}}})'\nonumber\\
\hspace{-1mm}&&\times\hspace{-1mm}
\Bigg(\prod^{d}_{j=0}\bar{A}_{\theta_{k\hspace{-0.51mm}-\hspace{-0.51mm}j}}z_{k\hspace{-0.51mm}-\hspace{-0.51mm}d}
\hspace{-1mm}+\hspace{-1mm}\sum^{d}_{i=1}\prod^{i-1}_{j=0}\bar{A}_{\theta_{k\hspace{-0.51mm}-\hspace{-0.51mm}j}}
\bar{B}_{\theta_{k\hspace{-0.51mm}-\hspace{-0.51mm}i}}u^{c}_{k\hspace{-0.51mm}
-\hspace{-0.51mm}d\hspace{-0.51mm}-\hspace{-0.51mm}i}\Bigg)\Bigg\}\nonumber\\
\hspace{-1mm}&=&\hspace{-1mm}
\mathbb{E}\Bigg\{z_{k}'Qz_{k}+(u^{c}_{k\hspace{-0.51mm}-\hspace{-0.51mm}d})'Ru^{c}_{k\hspace{-0.51mm}-\hspace{-0.51mm}d}
-(u^{c}_{k\hspace{-0.51mm}-\hspace{-0.51mm}d})'\Gamma_{\theta_{k\hspace{-0.51mm}
-\hspace{-0.51mm}d\hspace{-0.51mm}-\hspace{-0.51mm}1}}u^{c}_{k\hspace{-0.51mm}-\hspace{-0.51mm}d}\nonumber\\
\hspace{-1mm}&&-\hspace{-1mm}
(u^{c}_{k\hspace{-0.51mm}-\hspace{-0.51mm}d})'M^{0}_{\theta_{k\hspace{-0.51mm}-\hspace{-0.51mm}
d\hspace{-0.51mm}-\hspace{-0.51mm}1}}z_{k\hspace{-0.51mm}-\hspace{-0.51mm}d}
-(u^{c}_{k\hspace{-0.51mm}-\hspace{-0.51mm}d})'\sum^{d}_{i=1}M^{i}_{\theta_{k\hspace{-0.51mm}
-\hspace{-0.51mm}d\hspace{-0.51mm}-\hspace{-0.51mm}1}}
u^{c}_{k\hspace{-0.51mm}-\hspace{-0.51mm}d\hspace{-0.51mm}-\hspace{-0.51mm}i}\nonumber\\
\hspace{-1mm}&&-\hspace{-1mm}
z_{k\hspace{-0.51mm}-\hspace{-0.51mm}d}'(M^{0}_{\theta_{k\hspace{-0.51mm}-\hspace{-0.51mm}
d\hspace{-0.51mm}-\hspace{-0.51mm}1}})'u^{c}_{k\hspace{-0.51mm}-\hspace{-0.51mm}d}
\hspace{-1mm}-\hspace{-1mm}
z_{k\hspace{-0.51mm}-\hspace{-0.51mm}d}'
(M^{0}_{\theta_{k\hspace{-0.51mm}-\hspace{-0.51mm}d\hspace{-0.51mm}-\hspace{-0.51mm}1}})'
\Gamma^{-1}_{\theta_{k\hspace{-0.51mm}-\hspace{-0.51mm}d\hspace{-0.51mm}-\hspace{-0.51mm}1}}
M^{0}_{\theta_{k\hspace{-0.51mm}-\hspace{-0.51mm}d\hspace{-0.51mm}-\hspace{-0.51mm}1}}
z_{k\hspace{-0.51mm}-\hspace{-0.51mm}d}\nonumber\\
\hspace{-0.51mm}&&-\hspace{-0.51mm}z_{k\hspace{-0.51mm}-\hspace{-0.51mm}d}'\sum^{d}_{i=1}
(M^{0}_{\theta_{k\hspace{-0.51mm}-\hspace{-0.51mm}d\hspace{-0.51mm}-\hspace{-0.51mm}1}})'
\Gamma^{-1}_{\theta_{k\hspace{-0.51mm}-\hspace{-0.51mm}d\hspace{-0.51mm}-\hspace{-0.51mm}1}}
M^{i}_{\theta_{k\hspace{-0.51mm}-\hspace{-0.51mm}d\hspace{-0.51mm}-\hspace{-0.51mm}1}}
u^{c}_{k\hspace{-0.51mm}-\hspace{-0.51mm}d\hspace{-0.51mm}-\hspace{-0.51mm}i}\hspace{-1mm}-\hspace{-1mm}
\sum^{d}_{i=1}(u^{c}_{k\hspace{-0.51mm}-\hspace{-0.51mm}
d\hspace{-0.51mm}-\hspace{-0.51mm}i})'\nonumber\\
\hspace{-1mm}&&\times\hspace{-1mm}
(M^{i}_{\theta_{k\hspace{-0.51mm}-\hspace{-0.51mm}
d\hspace{-0.51mm}-\hspace{-0.51mm}1}})'u^{c}_{k\hspace{-0.51mm}-\hspace{-0.51mm}d}
\hspace{-0.51mm}-\hspace{-0.51mm}\sum^{d}_{i=1}(u^{c}_{k\hspace{-0.51mm}-\hspace{-0.51mm}
d\hspace{-0.51mm}-\hspace{-0.51mm}i})'(M^{i}_{\theta_{k\hspace{-0.51mm}-\hspace{-0.51mm}
d\hspace{-0.51mm}-\hspace{-0.51mm}1}})'\Gamma^{-1}_{\theta_{k\hspace{-0.51mm}-\hspace{-0.51mm}
d\hspace{-0.51mm}-\hspace{-0.51mm}1}}\nonumber\\
\hspace{-1mm}&&\times\hspace{-1mm}
\sum^{d}_{i=1}M^{i}_{\theta_{k\hspace{-0.51mm}-\hspace{-0.51mm}d\hspace{-0.51mm}
-\hspace{-0.51mm}1}}u^{c}_{k\hspace{-0.51mm}-\hspace{-0.51mm}d\hspace{-0.51mm}-\hspace{-0.51mm}i}
\Bigg\}\nonumber\\
\hspace{-1mm}&=&\hspace{-1mm}
\mathbb{E}\Bigg\{z_{k}'Qz_{k}+(u^{c}_{k\hspace{-0.51mm}-\hspace{-0.51mm}d})'R
u^{c}_{k\hspace{-0.51mm}-\hspace{-0.51mm}d}\Big(u^{c}_{k\hspace{-0.51mm}-\hspace{-0.51mm}d}
\hspace{-0.51mm}+\hspace{-0.51mm}\Gamma^{-1}_{\theta_{k\hspace{-0.51mm}-\hspace{-0.51mm}d
\hspace{-0.51mm}\hspace{-0.51mm}1}}M^{0}_{\theta_{k\hspace{-0.51mm}-\hspace{-0.51mm}d\hspace{-0.51mm}
-\hspace{-0.51mm}1}}
z_{k\hspace{-0.51mm}-\hspace{-0.51mm}d}\nonumber\\
\hspace{-1mm}&&+\hspace{-1mm}
\Gamma^{-1}_{\theta_{k\hspace{-0.51mm}-\hspace{-0.51mm}
d\hspace{-0.51mm}-\hspace{-0.51mm}1}}\sum^{d}_{i=1}M^{i}_{\theta_{k\hspace{-0.51mm}
-\hspace{-0.51mm}d\hspace{-0.51mm}-\hspace{-0.51mm}1}}
u^{c}_{k\hspace{-0.51mm}-\hspace{-0.51mm}d\hspace{-0.51mm}-\hspace{-0.51mm}i}\Big)'
\Gamma_{\theta_{k\hspace{-0.51mm}-\hspace{-0.51mm}d\hspace{-0.51mm}-\hspace{-0.51mm}1}}\Big(u^{c}_{k\hspace{-0.51mm}-\hspace{-0.51mm}d}\hspace{-0.51mm}+\hspace{-0.51mm}
\Gamma^{-1}_{\theta_{k\hspace{-0.51mm}-\hspace{-0.51mm}d\hspace{-0.51mm}-\hspace{-0.51mm}1}}\nonumber\\
\hspace{-1mm}&&\times\hspace{-1mm}
M^{0}_{\theta_{k\hspace{-0.51mm}-\hspace{-0.51mm}d\hspace{-0.51mm}-\hspace{-0.51mm}1}}
z_{k\hspace{-0.51mm}-\hspace{-0.51mm}d}
+\Gamma^{-1}_{\theta_{k\hspace{-0.51mm}-\hspace{-0.51mm}d\hspace{-0.51mm}-\hspace{-0.51mm}1}}
\sum^{d}_{i=1}M^{i}_{\theta_{k\hspace{-0.51mm}-\hspace{-0.51mm}d\hspace{-0.51mm}
-\hspace{-0.51mm}1}}
u^{c}_{k\hspace{-0.51mm}-\hspace{-0.51mm}d\hspace{-0.51mm}-\hspace{-0.51mm}i}\Big)\Bigg\}.
\label{f108}
\end{eqnarray}
\end{small}
Summing up from $k=d$ to $k=N$ on both sides of (\ref{f108}), and in view of $\Gamma_{\theta_{k-d-1}}>0$ for $k\geq d$,  the optimal controller and the optimal cost functional can be obtained as (\ref{f14}) and (\ref{f15}), respectively. The proof is completed.\\
\section{Proof of Theorem 2 \label{ap3}}

Proof.  ($\Longrightarrow$) \ \ The existence of solution to CAREs (\ref{f17})-(\ref{f21}) will be shown.  Consider the following delay-free MJLS:
\begin{eqnarray}
Y_{k+1}=C_{\theta_{k}}Y_{k}+Du^{c}_{k},\label{f36}
\end{eqnarray}
where
\begin{small}
\begin{eqnarray*}
Y_{k}=\left[
  \begin{array}{c}
z_{k}\\
u^{c}_{k-1}\\
\vdots\\
u^{c}_{k-d}
  \end{array}
\right],C_{\theta_{k}}=\left[
  \begin{array}{ccccc}
\bar{A}_{\theta_{k}}&0&\cdots&0&\bar{B}_{\theta_{k}}\\
0&0&\cdots&0&0\\
0&I&\cdots&0&0\\
\vdots&\vdots&\ddots&\vdots&\vdots\\
0&0&\cdots&I&0
  \end{array}
\right],
D=\left[
  \begin{array}{c}
0\\
I\\
0\\
\vdots\\
0
\end{array}
\right].
\end{eqnarray*}
\end{small}
Also, the cost functional over an infinite horizon is as follows
\begin{eqnarray}
\mathcal{J}=\sum_{k=0}^{\infty}\mathbb{E}[Y_{k}'\mathcal{Q}Y_{k}
+(u^{c}_{k})'Ru^{c}_{k}],\label{f38}
\end{eqnarray}
where $\mathcal{Q}=\left[
  \begin{array}{ccccc}
Q&  & &\\
 & 0 & &\\
 &  & \ddots & \\
 &  &  & 0
 \end{array}
\right].$
The corresponding cost functional over a finite horizon is
\begin{eqnarray}
\mathcal{J}_{N}\hspace{-1mm}=\hspace{-1mm}\mathbb{E}\Bigg\{\sum_{k=0}^{N}[Y_{k}'\mathcal{Q}Y_{k}
\hspace{-0.5mm}+\hspace{-0.5mm}(u^{c}_{k})'Ru^{c}_{k}]\hspace{-0.5mm}+\hspace{-0.5mm}
Y_{N\hspace{-0.5mm}+\hspace{-0.5mm}1}'\mathcal{P}_{N\hspace{-0.5mm}+\hspace{-0.5mm}1}
Y_{N\hspace{-0.5mm}+\hspace{-0.5mm}1}\Bigg\}.\label{f39}
\end{eqnarray}
By maximum principle, the following forward and backward difference equations can be given as
\begin{eqnarray}
\left\{
\begin{array}{lll}
0=\mathbb{E}_{k-1}[Ru^{c}_{k}+D'\xi_{k}],\\
\xi_{k-1}=\mathbb{E}_{k-1}[\mathcal{Q}Y_{k}+C_{\theta_{k}}\xi_{k}],\\
\xi_{N}=\mathcal{P}_{N+1}Y_{N+1}.
\end{array}
\right.\label{f40}
\end{eqnarray}
It follows from Theorem 1 that we have \\
(1) \ \  the following recursive sequence
\begin{small}
\begin{eqnarray}
\mathcal{P}^{(N)}_{N+1}\hspace{-1mm}&=&\hspace{-1mm}\bar{P}^{(N)}_{N+1},\label{f41}\\
\mathcal{P}^{(N)}_{\theta_{k\hspace{-0.51mm}-\hspace{-0.51mm}1}}\hspace{-1mm}&=&\hspace{-1mm}\mathbb{E}_{k\hspace{-0.51mm}-\hspace{-0.51mm}1}[\mathcal{Q}
\hspace{-1mm}+\hspace{-1mm}C_{\theta_{k}}'\mathcal{P}^{(N)}_{\theta_{k}}C_{\theta_{k}}
\hspace{-1mm}-\hspace{-1mm}(\mathcal{M}^{(N)}_{\theta_{k\hspace{-0.51mm}-\hspace{-0.51mm}1}})'
(\Upsilon^{(N)}_{\theta_{k\hspace{-0.51mm}-\hspace{-0.51mm}1}})^{-1}
\mathcal{M}^{(N)}_{\theta_{k-1}}],\label{f42}
\end{eqnarray}
\end{small}
in which
\begin{small}
\begin{eqnarray}
\mathcal{M}^{(N)}_{\theta_{k\hspace{-0.51mm}-\hspace{-0.51mm}1}}\hspace{-1mm}&=&\hspace{-1mm}\mathbb{E}_{k\hspace{-0.51mm}-\hspace{-0.51mm}1}[
D'\mathcal{P}^{(N)}_{\theta_{k}}C_{\theta_{k}}]\nonumber\\
\hspace{-1mm}&=&\hspace{-1mm}\mathbb{E}_{k\hspace{-0.51mm}-\hspace{-0.51mm}1}[(\mathcal{P}^{(N)}_{\theta_{k}}(2,1)
\bar{A}_{\theta_{k}}, \mathcal{P}^{(N)}_{\theta_{k}}(2,3),\hspace{-1mm}\cdots,\hspace{-1mm}
\mathcal{P}^{(N)}_{\theta_{k}}(2,d+1),\nonumber\\
\hspace{-1mm}&&\hspace{-1mm}\mathcal{P}^{(N)}_{\theta_{k}}(2,1)\bar{B}_{\theta_{k}})],\label{f44}\\
\Upsilon^{(N)}_{\theta_{k-1}}&=&\mathbb{E}_{k\hspace{-0.51mm}-\hspace{-0.51mm}1}[R
\hspace{-0.51mm}+\hspace{-0.51mm}D'\mathcal{P}^{(N)}_{\theta_{k\hspace{-0.51mm}-\hspace{-0.51mm}1}}D]
\hspace{-1mm}=\hspace{-1mm}\mathbb{E}_{k\hspace{-0.51mm}-\hspace{-0.51mm}1}[R\hspace{-0.51mm}+\hspace{-0.51mm}\mathcal{P}^{(N)}_{\theta_{k\hspace{-0.51mm}-\hspace{-0.51mm}1}}(2,2)];\label{f044}
\end{eqnarray}
\end{small}
(2) \ the costate
\begin{eqnarray}
\xi_{k-1}=\mathcal{P}^{(N)}_{\theta_{k-1}}Y_{k};\label{f450}
\end{eqnarray}
(3) \ the optimal control
\begin{small}
\begin{eqnarray}
u^{c}_{k}\hspace{-1mm}&=&\hspace{-1mm}-(\Upsilon^{(N)}_{\theta_{k\hspace{-0.51mm}-\hspace{-0.51mm}1}})^{-1}
\mathcal{M}^{(N)}_{\theta_{k\hspace{-0.51mm}-\hspace{-0.51mm}1}}Y_{k}\nonumber\\
\hspace{-1mm}&=&\hspace{-1mm}-(R+\mathcal{P}^{(N)}_{\theta_{k\hspace{-0.51mm}-\hspace{-0.51mm}1}}(2,2))^{-1}\big[
(\mathcal{P}^{(N)}_{\theta_{k}}(2,1)\bar{A}_{\theta_{k}})z_{k}\hspace{-0.51mm}+ \hspace{-0.51mm} \mathcal{P}^{(N)}_{\theta_{k}}(2,3)\nonumber\\
\hspace{-1mm}&&\times\hspace{-1mm}u^{c}_{k\hspace{-0.51mm}-\hspace{-0.51mm}1}
\hspace{-1mm}+\cdots+\hspace{-1mm}\mathcal{P}^{(N)}_{\theta_{k}}(2,d\hspace{-0.51mm}+\hspace{-0.51mm}1)u^{c}_{k\hspace{-0.51mm}-\hspace{-0.51mm}d\hspace{-0.51mm}+\hspace{-0.51mm}1}
\hspace{-1mm}+\hspace{-1mm}\mathcal{P}^{(N)}_{\theta_{k}}(2,1)\bar{B}_{\theta_{k}}u^{c}_{k\hspace{-0.51mm}-\hspace{-0.51mm}d}
\big]\nonumber\\
\hspace{-1mm}&=&\hspace{-1mm}-\Gamma^{-1}_{\theta_{k\hspace{-0.51mm}-\hspace{-0.51mm}1}}(k,N)M^{0}_{\theta_{k\hspace{-0.51mm}-\hspace{-0.51mm}1}}(k,N)z_{k}
\hspace{-0.51mm}-\hspace{-0.51mm}\Gamma^{-1}_{\theta_{k-1}}(k,N)\nonumber\\
\hspace{-1mm}&&\times\hspace{-1mm}
\sum^{d}_{i=1}M^{i}_{\theta_{k\hspace{-0.51mm}-\hspace{-0.51mm}1}}(k,N)u^{c}_{k\hspace{-0.51mm}-\hspace{-0.51mm}i},\label{f45}
\end{eqnarray}
\end{small}
where $\mathcal{P}^{(N)}_{\theta_{k}}(i,j)$ represents block matrix with suitable dimension.\\
The following relationship can be obtained in view of (\ref{f14}) and (\ref{f45}).
\begin{small}\begin{eqnarray}
\left\{
\begin{array}{lll}
&&(R+\mathcal{P}^{(N)}_{\theta_{k-1}}(2,2))^{-1}\mathbb{E}_{k-1}[
 \mathcal{P}^{(N)}_{\theta_{k}}(2,1)\bar{A}_{\theta_{k}}]\nonumber\\
 &=&
\Gamma^{-1}_{\theta_{k-1}}(k,N)M^{0}_{\theta_{k-1}}(k,N)\\
&&(R+\mathcal{P}^{(N)}_{\theta_{k-1}}(2,2))^{-1}\mathbb{E}_{k-1}[\mathcal{P}^{(N)}_{\theta_{k}}
(2,3)]\\
&=&
\Gamma^{-1}_{\theta_{k-1}}(k,N)M^{1}_{\theta_{k-1}}(k,N)\\
\ \ \ &\vdots&\\
&&(R+\mathcal{P}^{(N)}_{\theta_{k-1}}(2,2))^{-1}\mathbb{E}_{k-1}[\mathcal{P}^{(N)}_{\theta_{k}}
(2,d+1)]\\
&=&
\Gamma^{-1}_{\theta_{k-1}}(k,N)M^{d-1}_{\theta_{k-1}}(k,N)\\
&&(R+\mathcal{P}^{(N)}_{\theta_{k-1}}(2,2))^{-1}\mathbb{E}_{k-1}[\mathcal{P}^{(N)}_{\theta_{k}}(2,1)
\bar{B}_{\theta_{k}}]\\
&=&\Gamma^{-1}_{\theta_{k-1}}(k,N)M^{d}_{\theta_{k-1}}(k,N).\label{f46}
\end{array}
\right.
\end{eqnarray}\end{small}
The convergence of $\mathcal{P}^{(N)}_{\theta_{k-1}}$ can be obtained  in a similar manner with \cite{21}. On this basis, from  (\ref{f46}), $\mathcal{P}^{(N)}_{\theta_{k-1}}$ $\Gamma^{-1}_{\theta_{k-1}}(k,N)M^{j}_{\theta_{k-1}}(k,N), j=0,1,\cdots, d$ are convergent.
Let $\xi_{k-1}=\left[
  \begin{array}{cccc}
\xi^{0}_{k-1}&\xi^{1}_{k-1}&\vdots&\xi^{d}_{k-1}
  \end{array}
\right]'$, and from (\ref{f450}) we know
\begin{eqnarray}
\xi^{0}_{k-1}&=&\mathcal{P}^{(N)}_{\theta_{k-1}}(1,1)z_{k}
+\mathcal{P}^{(N)}_{\theta_{k-1}}(1,2)u^{c}_{k-1}
+\cdots\nonumber\\
&&+\mathcal{P}^{(N)}_{\theta_{k-1}}(1,d+1)u^{c}_{k-d}.\label{f47}
\end{eqnarray}
Further, from (\ref{f40}), we have
\begin{eqnarray}
\xi^{0}_{k-1}=\mathbb{E}_{k-1}[Qz_{k}
+\bar{A}_{\theta_{k}}\xi^{0}_{k}],\label{f48}
\end{eqnarray}
comparing with (\ref{f5}) and (\ref{f48}), it is obvious that, if $\xi^{0}_{N}=\bar{P}_{N+1}$, the following relationship holds
\begin{eqnarray}
\xi^{0}_{k-1}=\lambda_{k-1}. \label{f49}
\end{eqnarray}
Consider (\ref{f14}), (\ref{f16}), (\ref{f47}) and (\ref{f49}), and we can find the following relationship using a direct calculation
\begin{eqnarray}
\mathcal{P}^{(N)}_{\theta_{k-1}}(1,1)
=\bar{P}_{\theta_{k-1}}.\label{f50}
\end{eqnarray}
Therefore, we have
\begin{eqnarray}
\lim\limits_{N\rightarrow\infty}\bar{P}_{\theta_{k-1}}(k,N)\triangleq \bar{P}_{l_{d-1}},\label{f500}
\end{eqnarray}
where $\theta_{k-1}=l_{d-1}, k\geq d, l_{d-1}=0,1.$
In view of (\ref{f8}), we know that $(M^{0}_{\theta_{k-1}})'\Gamma^{-1}_{\theta_{k-1}}M^{0}_{\theta_{k-1}}$ is convergent.\\
On this basis, it is esay to verify that $\Gamma_{\theta_{k-1}}(k,N)$, $M^{i}_{\theta_{k-1}}(k,N),i=0,\cdots, d$,  $F^{d-j+1}_{\theta_{k-j-1}}(k-j,N)$ and $\tilde{S}^{j}_{\theta_{k-1}}(k,N),j=1,\cdots, d$ are also convergent and  CAREs (\ref{f17})-(\ref{f21}) have a solution.\\
(2) \ Now, we will prove inequality (\ref{f22}).\\
Let $\hat{J}_{d}(m)$ represent cost functional (\ref{f3}) which starts at $d$ and ends at $m, m\geq N$. In view of Lemma 1, we can obtain the optimal cost value $\hat{J}^{\ast}_{d}(m)$ as follows
\begin{eqnarray}
\hspace{-1mm}&&\hspace{-1mm}\hat{J}^{\ast}_{d}(m)\nonumber\\
\hspace{-1mm}&=&\hspace{-1mm}\mathbb{E}\Bigg\{z_{d}'\bar{P}_{\theta_{d\hspace{-0.5mm}-\hspace{-0.5mm}1}}(d,m)z_{d}
\hspace{-1mm}-\hspace{-1mm}z_{d}'\sum_{s=0}^{d-1}[F^{s\hspace{-0.5mm}+\hspace{-0.5mm}1}_{\theta_{d\hspace{-0.5mm}-\hspace{-0.5mm}1}}(d,m)
\Gamma^{-1}_{\theta_{s\hspace{-0.5mm}-\hspace{-0.5mm}1}}(d+s,m)\nonumber\\
\hspace{-1mm}&&\times\hspace{-1mm}
M^{0}_{\theta_{s\hspace{-0.5mm}-\hspace{-0.5mm}1}}(d\hspace{-0.5mm}+\hspace{-0.5mm}s)z_{s}]
\hspace{-1mm}-\hspace{-1mm}z_{d}'\sum_{s=0}^{d-1}\Bigg[F^{s\hspace{-0.5mm}+\hspace{-0.5mm}1}_{\theta_{d\hspace{-0.5mm}-\hspace{-0.5mm}1}}(d,m)
\Gamma^{-1}_{\theta_{s\hspace{-0.5mm}-\hspace{-0.5mm}1}}(d\hspace{-0.5mm}+\hspace{-0.5mm}s,m)\nonumber\\
\hspace{-1mm}&&\times\hspace{-1mm}
\sum_{i=s+1}^{d}M^{i}_{\theta_{s\hspace{-0.5mm}-\hspace{-0.5mm}1}}(d\hspace{-0.5mm}
+\hspace{-0.5mm}s)u^{c}_{s\hspace{-0.5mm}-\hspace{-0.5mm}i}\Bigg]\Bigg\}\nonumber
\end{eqnarray}
\begin{eqnarray}
\hspace{-1mm}&=&\hspace{-1mm}\mathbb{E}\Bigg\{z_{d}'\Bigg[\bar{P}_{\theta_{d\hspace{-0.5mm}-\hspace{-0.5mm}1}}(d,m)
\hspace{-1mm}-\hspace{-1mm}\sum_{s=0}^{d-1}\Big[F^{s\hspace{-0.5mm}+\hspace{-0.5mm}1}_{\theta_{d\hspace{-0.5mm}-\hspace{-0.5mm}1}}(d,m)
\Gamma^{-1}_{\theta_{s\hspace{-0.5mm}-\hspace{-0.5mm}1}}(d\hspace{-0.5mm}+\hspace{-0.5mm}s,m)\nonumber\\
\hspace{-1mm}&&\times\hspace{-1mm}
(F^{s\hspace{-0.5mm}+\hspace{-0.5mm}1}_{\theta_{d\hspace{-0.5mm}-\hspace{-0.5mm}1}}(d,m))'\Big]\Bigg]z_{d}\Bigg\}\nonumber\\
\hspace{-1mm}&=&\hspace{-1mm}z_{d}'\Bigg[\bar{P}_{\theta_{d\hspace{-0.5mm}-\hspace{-0.5mm}1}}(d,m)
\hspace{-1mm}-\hspace{-1mm}\sum_{s=0}^{d-1}\Big[F^{s\hspace{-0.5mm}+\hspace{-0.5mm}1}_{\theta_{d\hspace{-0.5mm}-\hspace{-0.5mm}1}}(d,m)
\Gamma^{-1}_{\theta_{s\hspace{-0.5mm}-\hspace{-0.5mm}1}}(d\hspace{-0.5mm}+\hspace{-0.5mm}s,m)
\nonumber\\
\hspace{-1mm}&&\times\hspace{-1mm}
(F^{s\hspace{-0.5mm}+\hspace{-0.5mm}1}_{\theta_{d\hspace{-0.5mm}-\hspace{-0.5mm}1}}(d,m))'\Big]\Bigg]z_{d} \nonumber\\
 \hspace{-1mm}&\geq&\hspace{-1mm} 0. \label{f0023}
\end{eqnarray}
Since $z_{d}$ is arbitrary, we have
\begin{eqnarray}
\hspace{-1mm}&&\hspace{-0.5mm}\bar{P}_{\theta_{d\hspace{-0.5mm}-\hspace{-0.5mm}1}}(d,m)
\hspace{-1mm}-\hspace{-1mm}\sum_{s=0}^{d-1}\Big[F^{s\hspace{-0.5mm}+\hspace{-0.5mm}1}_{\theta_{d\hspace{-0.5mm}-\hspace{-0.5mm}1}}(d,m)
\Gamma^{-1}_{\theta_{s\hspace{-0.5mm}-\hspace{-0.5mm}1}}(d\hspace{-0.5mm}+\hspace{-0.5mm}s,m)
\nonumber\\
\hspace{-1mm}&&\times\hspace{-1mm}
(F^{s\hspace{-0.5mm}+\hspace{-0.5mm}1}_{\theta_{d\hspace{-0.5mm}-\hspace{-0.5mm}1}}(d,m))'\Big]\geq 0.\label{f0021}
\end{eqnarray}
For $k\geq d$, let $m=N-k+d$. In view of the time-variance, it yields that
\begin{eqnarray}
\hspace{-1mm}&&\hspace{-1mm}\bar{P}_{\theta_{k\hspace{-0.5mm}-\hspace{-0.5mm}1}}(k,N)
-\sum_{s=0}^{d-1}\Big[F^{s\hspace{-0.5mm}+\hspace{-0.5mm}1}_{\theta_{k\hspace{-0.5mm}-\hspace{-0.5mm}1}}(k,N)\Gamma^{-1}_{\theta_{k\hspace{-0.5mm}-\hspace{-0.5mm}d\hspace{-0.5mm}+\hspace{-0.5mm}s\hspace{-0.5mm}-\hspace{-0.5mm}1}}(k\hspace{-0.5mm}+\hspace{-0.5mm}s,N)
\nonumber\\
\hspace{-1mm}&&\times\hspace{-1mm}(F^{s\hspace{-0.5mm}+\hspace{-0.5mm}1}_{\theta_{k\hspace{-0.5mm}-\hspace{-0.5mm}1}}(k,N))'\Big]
\geq 0.\label{f0022}
\end{eqnarray}
By virtue of the convergence, it is easy to derive (\ref{f22}). \\
From Lemma 3 in \cite{18}, we can find an integer $G$ such that
\begin{eqnarray}
\hspace{-1mm}&&\hspace{-1mm}\bar{P}_{\theta_{d}}(d,G)
-\sum_{s=0}^{d-1}[F^{s+1}_{\theta_{d}}(d,G)\Gamma^{-1}_{\theta_{s-1}}(d+s,G)
\nonumber\\
\hspace{-1mm}&&\times\hspace{-1mm}(F^{s+1}_{\theta_{d}}(d,G))']> 0. \label{f0900}
\end{eqnarray}
Also, the monotonicity with respect ro $N$ of (\ref{f0022}) deduces that
\begin{eqnarray*}
\hspace{-1mm}&&\hspace{-1mm}\bar{P}_{l_{d}}
\hspace{-1mm}-\hspace{-1mm}\sum_{s=0}^{d-1}[F^{s+1}_{l_{d}}\Gamma^{-1}_{l_{s-1}}
(F^{s+1}_{l_{d}})']\\
\hspace{-1mm}&=&\hspace{-1mm}\underset{N\rightarrow\infty}{\lim}\Bigg[\bar{P}_{\theta_{d}}(d,N)
\hspace{-1mm}-\hspace{-1mm}\sum_{s=0}^{d-1}\Big(F^{s+1}_{\theta_{d}}(d,N)\Gamma^{-1}_{\theta_{s-1}}(d+s,N)
\nonumber\\
\hspace{-1mm}&&\times\hspace{-1mm}(F^{s+1}_{\theta_{d}}(d,N))'\Big)\Bigg]\\
\hspace{-1mm}&\geq&\hspace{-1mm} \bar{P}_{\theta_{d}}(d,G)
\hspace{-1mm}-\hspace{-1mm}\sum_{s=0}^{d-1}\Big[F^{s+1}_{\theta_{d}}(d,G)\Gamma^{-1}_{\theta_{s-1}}(d+s,G)
\nonumber\\
\hspace{-1mm}&&\times\hspace{-1mm}(F^{s+1}_{\theta_{d}}(d,G))'\Big]>0.
\end{eqnarray*}
The proof is completed.\\
($\Longleftarrow$) \ \ The mean-square stabilization of system (\ref{f1}) will be illustrated. Define
\begin{eqnarray}
\mathcal{L}(k)\hspace{-1mm}&=&\hspace{-1mm}\mathbb{E}\Bigg\{\hspace{-0.5mm}z_{k}'\bar{P}_{\theta_{k\hspace{-0.5mm}-\hspace{-0.5mm}1}}
z_{k}\hspace{-1mm}-\hspace{-1mm}z_{k}'\sum^{d-1}_{s=0}\Big(F^{s\hspace{-0.5mm}+\hspace{-0.5mm}1}_{\theta_{k\hspace{-0.5mm}-\hspace{-0.5mm}1}}
\Gamma^{-1}_{\theta_{k\hspace{-0.5mm}-\hspace{-0.5mm}d\hspace{-0.5mm}-\hspace{-0.5mm}1\hspace{-0.5mm}+\hspace{-0.5mm}s}}M^{0}_{\theta_{k\hspace{-0.5mm}-\hspace{-0.5mm}d\hspace{-0.5mm}-\hspace{-0.5mm}1\hspace{-0.5mm}+\hspace{-0.5mm}s}}z_{k\hspace{-0.5mm}-\hspace{-0.5mm}d\hspace{-0.5mm}+\hspace{-0.5mm}s}\Big)\nonumber\\
\hspace{-1mm}&&-\hspace{-1mm}z_{k}'\sum^{d-1}_{s=0}\Big(F^{s\hspace{-0.5mm}+\hspace{-0.5mm}1}_{\theta_{k\hspace{-0.5mm}-\hspace{-0.5mm}1}}
\Gamma^{-1}_{\theta_{k\hspace{-0.5mm}-\hspace{-0.5mm}d\hspace{-0.5mm}
-\hspace{-0.5mm}1\hspace{-0.5mm}+\hspace{-0.5mm}s}}\sum^{d}_{i=s+1}
M^{i}_{\theta_{k\hspace{-0.5mm}-\hspace{-0.5mm}d\hspace{-0.5mm}-\hspace{-0.5mm}1\hspace{-0.5mm}
+\hspace{-0.5mm}s}}u^{c}_{k\hspace{-0.5mm}-\hspace{-0.5mm}d\hspace{-0.5mm}
-\hspace{-0.5mm}i\hspace{-0.5mm}+\hspace{-0.5mm}s}\Big)\hspace{-0.5mm}\Bigg\}.\label{f58}
\end{eqnarray}
In view of (\ref{f105}),  we have
\begin{eqnarray}
\hspace{-1mm}&&\hspace{-1mm}\mathcal{L}(k)-\mathcal{L}(k+1)\nonumber\\
\hspace{-1mm}&=&\hspace{-1mm}\mathbb{E}\Bigg\{z_{k}'Qz_{k}\hspace{-1mm}+\hspace{-1mm}(u^{c}_{k\hspace{-0.5mm}-\hspace{-0.5mm}d})'Ru^{c}_{k\hspace{-0.5mm}-\hspace{-0.5mm}d}-\Bigg(u^{c}_{k\hspace{-0.5mm}-\hspace{-0.5mm}d}
\hspace{-1mm}+\hspace{-1mm}\Gamma^{-1}_{i}M^{0}_{i}z_{k\hspace{-0.5mm}-\hspace{-0.5mm}d}
\hspace{-1mm}+\hspace{-1mm}\Gamma^{-1}_{i}\nonumber\\
\hspace{-1mm}&&\times\hspace{-1mm}\sum^{d}_{s=1}
(M^{s}_{i}u^{c}_{k\hspace{-0.5mm}-\hspace{-0.5mm}d\hspace{-0.5mm}-\hspace{-0.5mm}s})\Bigg)'\Gamma_{i}
\Bigg(u^{c}_{k\hspace{-0.5mm}-\hspace{-0.5mm}d}\hspace{-1mm}+\hspace{-1mm}\
\Gamma^{-1}_{i}M^{0}_{i}z_{k\hspace{-0.5mm}-\hspace{-0.5mm}d}\nonumber\\
\hspace{-1mm}&&+\hspace{-1mm}
\Gamma^{-1}_{i}\sum^{d}_{s=1}
(M^{s}_{i}
u^{c}_{k\hspace{-0.5mm}-\hspace{-0.5mm}d\hspace{-0.5mm}-\hspace{-0.5mm}s})\Bigg)\Bigg\}\label{f059}\\
\hspace{-1mm}&=&\hspace{-1mm}\mathbb{E}\Big[z_{k}'Qz_{k}
\hspace{-1mm}+\hspace{-1mm}(u^{c}_{k\hspace{-0.5mm}-\hspace{-0.5mm}d})'R
u^{c}_{k\hspace{-0.5mm}-\hspace{-0.5mm}d}\Big]\geq0, \ \ k\geq d,\label{f59}
\end{eqnarray}
with $u^{c}_{k-d}=-\Gamma^{-1}_{i}\Bigg(M^{0}_{i}z_{k-d}
+\sum^{d}_{s=1}M^{s}_{i}u^{c}_{k-d-s}\Bigg)$. Therefore, $\mathcal{L}(k)$  decreases with respect to $k$. From (\ref{e01})-(\ref{e0}), $\mathcal{L}(k)$ can be expressed as
\begin{eqnarray}
\mathcal{L}(k)\hspace{-1mm}&=&\hspace{-1mm}\mathbb{E}\Bigg\{z_{k}'\Bigg[\bar{P}_{\theta_{k-1}}
\hspace{-1mm}-\hspace{-1mm}\sum_{s=0}^{d-1}\Big(F^{s+1}_{\theta_{k-1}}\Gamma^{-1}_{\theta_{k-d+s-1}}
(F^{s+1}_{\theta_{k-1}})'\Big)\Bigg]z_{k}\nonumber\\
\hspace{-1mm}&&+\hspace{-1mm}\sum_{s=0}^{d-1}\Big(F^{s+1}_{\theta_{k-1}}z_{k}
\hspace{-1mm}-\hspace{-1mm}\mathbb{E}[F^{s+1}_{\theta_{k-1}}z_{k}
|\mathcal{F}_{k-s-1}]\Big)'\Gamma^{-1}_{\theta_{k-s-1}}\nonumber\\
\hspace{-1mm}&&\times\hspace{-1mm} \Big(F^{s+1}_{\theta_{k-1}}z_{k}-\mathbb{E}[F^{s+1}_{\theta_{k-1}}z_{k}
|\mathcal{F}_{k-s-1}]\Big)\Bigg\}\nonumber\\
\hspace{-1mm}&\geq&\hspace{-1mm}\mathbb{E}\Bigg\{z_{k}'\Bigg[\bar{P}_{\theta_{k-1}}
\hspace{-1mm}-\hspace{-1mm}\sum_{s=0}^{d-1}[F^{s+1}_{\theta_{k-1}}\Gamma^{-1}_{\theta_{k-d+s-1}}
(F^{s+1}_{\theta_{k-1}})']\Bigg]z_{k}\bigg\}\nonumber\\
\hspace{-1mm}&\geq&\hspace{-1mm}0, \ \ k\geq d, \label{f60}
\end{eqnarray}
i.e., $\mathcal{L}(k)$ is bounded. Therefore, $\mathcal{L}(k)$ is convergent.

For any $l\geq 0$, summing up from $k=l+d$ to $k=l+N$ on both sides of (\ref{f59}), when $l\rightarrow\infty$, we can derive that
\begin{eqnarray}
&&\lim\limits_{l\rightarrow\infty}\sum_{k=l+d}^{l+N}\mathbb{E}[z_{k}'Qz_{k}
+(u^{c}_{k-d})'Ru^{c}_{k-d}]\nonumber\\
&=&\lim\limits_{l\rightarrow\infty}[\mathcal{L}(l+d)-\mathcal{L}(l+N+1)]=0.\label{f61}
\end{eqnarray}
Recall that
\begin{eqnarray*}
&&\sum_{k=d}^{N}\mathbb{E}[z_{k}'Qz_{k}
+(u^{c}_{k-d})'Ru^{c}_{k-d}]\\
&\geq & \mathbb{E}\Bigg\{z_{d}'\bar{P}_{\theta_{d-1}}z_{d}-z_{d}'\sum^{d-1}_{s=0}\left(F^{s+1}_{\theta_{d-1}}
\Gamma^{-1}_{\theta_{s-1}}M^{0}_{\theta_{s-1}}z_{s}\right)\nonumber\\
&&-z_{d}'\sum^{d}_{s=0}\left(F^{s+1}_{\theta_{d-1}}
\Gamma^{-1}_{\theta_{s-1}}\sum^{d}_{i=s+1}M^{i}_{\theta_{s-1}}u^{c}_{s-i}\right)\Bigg\}.
\end{eqnarray*}
Therefore, the following relationship can be deduced that
\begin{eqnarray*}
\hspace{-1mm}&&\hspace{-1mm}\sum_{k=l+d}^{l+N}\mathbb{E}\Big[z_{k}'Qz_{k}
\hspace{-1mm}+\hspace{-1mm}(u^{c}_{k\hspace{-0.5mm}-\hspace{-0.5mm}d})'Ru^{c}_{k\hspace{-0.5mm}-\hspace{-0.5mm}d}\Big]\\
\hspace{-1mm}&\geq &\hspace{-1mm} \mathbb{E}\Bigg\{z_{l\hspace{-0.5mm}+\hspace{-0.5mm}d}'\Bigg[\bar{P}_{\theta_{l\hspace{-0.5mm}+\hspace{-0.5mm}d\hspace{-0.5mm}-\hspace{-0.5mm}1}}(l\hspace{-0.5mm}+\hspace{-0.5mm}d,l\hspace{-0.5mm}+\hspace{-0.5mm}N)
\hspace{-1mm}-\hspace{-1mm}\sum_{s=0}^{d-1}[F^{s\hspace{-0.5mm}
+\hspace{-0.5mm}1}_{\theta_{l\hspace{-0.5mm}+\hspace{-0.5mm}d\hspace{-0.5mm}
-\hspace{-0.5mm}1}}(l\hspace{-0.5mm}+\hspace{-0.5mm}d,l\hspace{-0.5mm}+\hspace{-0.5mm}N)\\
\hspace{-1mm}&&\times\hspace{-1mm}
\Gamma^{-1}_{\theta_{l\hspace{-0.5mm}+\hspace{-0.5mm}s\hspace{-0.5mm}-\hspace{-0.5mm}1}}
(l\hspace{-0.5mm}+\hspace{-0.5mm}s,l\hspace{-0.5mm}+\hspace{-0.5mm}N) (F^{s\hspace{-0.5mm}+\hspace{-0.5mm}1}_{\theta_{l\hspace{-0.5mm}+\hspace{-0.5mm}d\hspace{-0.5mm}-\hspace{-0.5mm}1}}(l\hspace{-0.5mm}+\hspace{-0.5mm}d,l\hspace{-0.5mm}+\hspace{-0.5mm}N))']\Bigg]z_{l\hspace{-0.5mm}+\hspace{-0.5mm}d}\Bigg\}\\
\hspace{-1mm}&=&\hspace{-1mm}\mathbb{E}\Bigg\{z_{l\hspace{-0.5mm}+\hspace{-0.5mm}d}'\Bigg[\bar{P}_{\theta_{d\hspace{-0.5mm}-\hspace{-0.5mm}1}}(d,N)
\hspace{-1mm}-\hspace{-1mm}\sum_{s=0}^{d-1}[F^{s\hspace{-0.5mm}+\hspace{-0.5mm}1}_{\theta_{d\hspace{-0.5mm}-\hspace{-0.5mm}1}}(d,N)
\Gamma^{-1}_{\theta_{s\hspace{-0.5mm}-\hspace{-0.5mm}1}}(s,N)\\
\hspace{-1mm}&&\times\hspace{-1mm}
(F^{s\hspace{-0.5mm}+\hspace{-0.5mm}1}_{\theta_{d\hspace{-0.5mm}-\hspace{-0.5mm}1}}(d,N))']\Bigg]z_{l\hspace{-0.5mm}+\hspace{-0.5mm}d}\Bigg\}\\
\hspace{-1mm}&\geq&\hspace{-1mm} 0.
\end{eqnarray*}
Using (\ref{f61}), we have
\begin{eqnarray}
\hspace{-1mm}&&\hspace{-1mm}\lim\limits_{l\rightarrow\infty}\mathbb{E}\Bigg\{z_{l\hspace{-0.5mm}+\hspace{-0.5mm}d}'
\Bigg[\bar{P}_{\theta_{d\hspace{-0.5mm}-\hspace{-0.5mm}1}}(d,N)
\hspace{-1mm}-\hspace{-1mm}\sum_{s=0}^{d-1}[F^{s\hspace{-0.5mm}+\hspace{-0.5mm}1}_{\theta_{d
\hspace{-0.5mm}-\hspace{-0.5mm}1}}(d,N)\nonumber\\
\hspace{-1mm}&&\times\hspace{-1mm}
\Gamma^{-1}_{\theta_{s\hspace{-0.5mm}-\hspace{-0.5mm}1}}(s,N)
(F^{s\hspace{-0.5mm}+\hspace{-0.5mm}1}_{\theta_{d\hspace{-0.5mm}-\hspace{-0.5mm}1}}(d,N))']\Bigg]z_{l\hspace{-0.5mm}+\hspace{-0.5mm}d}\Bigg\}\nonumber\\
\hspace{-1mm}&=&\hspace{-1mm}0, \ \ \forall N\geq d.\label{f62}
\end{eqnarray}
From (\ref{f0900}) and (\ref{f62}), it deduces that
$\lim\limits_{l\rightarrow\infty}\mathbb{E}[z_{l+d}'z_{l+d}]=0$.
Therefore, system (\ref{f3}) is mean-square stabilizable with controller (\ref{f25}).

Finally, the optimal cost functional will be calculated.

Summing up from $k=0$ to $k=N$ on both sides of (\ref{f059}), it yields that
\begin{eqnarray}
\hspace{-1mm}&&\hspace{-1mm}\mathbb{E}\Bigg\{\sum_{k=0}^{N}[z_{k}'Qz_{k}
\hspace{-1mm}+\hspace{-1mm}(u^{c}_{k-d})'Ru^{c}_{k-d}]\Bigg\}\nonumber\\
\hspace{-1mm}&=&\hspace{-1mm}\mathcal{L}(0)\hspace{-1mm}-\hspace{-1mm}\mathcal{L}(N+1)
\hspace{-1mm}+\hspace{-1mm}\sum_{k=0}^{N}\mathbb{E}\bigg\{\Big(u^{c}_{k-d}
\hspace{-1mm}+\hspace{-1mm}\Gamma^{-1}_{i}M^{0}_{i}z_{k-d}
\nonumber\\
\hspace{-1mm}&&+\hspace{-1mm}\Gamma^{-1}_{i}\sum^{d}_{s=1}(M^{s}_{i}u^{c}_{k-d-s})\Big)' \Gamma_{i}\Big(u^{c}_{k-d}\hspace{-1mm}+\hspace{-1mm}\Gamma^{-1}_{i}M^{0}_{i}z_{k-d}\nonumber\\
\hspace{-1mm}&&+\hspace{-1mm}
\Gamma^{-1}_{i}\sum^{d}_{s=1}(M^{s}_{i}u^{c}_{k-d-s})\Big)\bigg\}.\label{f63}
\end{eqnarray}
Since $0\leq\mathcal{L}(k)\leq \mathbb{E}[z_{k}\bar{P}_{\theta_{k-1}}z_{k}]$
 and the system (\ref{f5}) is stabilized in the mean-square sense, we have $\lim\limits_{k\rightarrow\infty}\mathbb{E}[z_{k}\bar{P}_{\theta_{k-1}}z_{k}]=0$, i.e.,
  $\lim\limits_{k\rightarrow\infty}\mathcal{L}(k)=0$.

Let $N\rightarrow\infty$ on both sides of (\ref{f63}), then
\begin{eqnarray}
J\hspace{-1mm}&=&\hspace{-1mm}\mathcal{L}(0)\hspace{-1mm}+\hspace{-1mm}\sum_{k=d}^{\infty}
\mathbb{E}\Bigg\{\Big(u^{c}_{k\hspace{-0.5mm}-\hspace{-0.5mm}d}
\hspace{-1mm}+\hspace{-1mm}\Gamma^{-1}_{i}M^{0}_{i}z_{k\hspace{-0.5mm}-\hspace{-0.5mm}d}
\hspace{-1mm}+\hspace{-1mm}\Gamma^{-1}_{i}\sum^{d}_{s=1}(M^{s}_{i}u^{c}_{k\hspace{-0.5mm}-\hspace{-0.5mm}d\hspace{-0.5mm}-\hspace{-0.5mm}s})\Big)'\nonumber\\
\hspace{-1mm}&&\times \hspace{-1mm}\Gamma_{i}\Big(u^{c}_{k\hspace{-0.5mm}-\hspace{-0.5mm}d}\hspace{-0.5mm}+\hspace{-0.5mm}\Gamma^{-1}_{i}M^{0}_{i}z_{k\hspace{-0.5mm}-\hspace{-0.5mm}d}
\hspace{-1mm}+\hspace{-1mm}\Gamma^{-1}_{i}\sum^{d}_{s=1}(M^{s}_{i}u^{c}_{k\hspace{-0.5mm}-\hspace{-0.5mm}d\hspace{-0.5mm}-\hspace{-0.5mm}s})\Big)\bigg\}\nonumber\\
\hspace{-1mm}&&+\hspace{-1mm}\sum_{k=0}^{d-1}\mathbb{E}
\bigg[\Big(u^{c}_{k-d}\hspace{-1mm}+\hspace{-1mm}\Gamma^{-1}_{i}M^{0}_{i}z_{k\hspace{-0.5mm}-\hspace{-0.5mm}d}
\hspace{-1mm}+\hspace{-1mm}\Gamma^{-1}_{i}\sum^{d}_{s=1}(M^{s}_{i}u^{c}_{k\hspace{-0.5mm}-\hspace{-0.5mm}d\hspace{-0.5mm}-\hspace{-0.5mm}s})\Big)'\nonumber\\
\hspace{-1mm}&&\times\hspace{-1mm} \Gamma_{i}\Big(u^{c}_{k\hspace{-0.5mm}-\hspace{-0.5mm}d}\hspace{-1mm}+\hspace{-1mm}\Gamma^{-1}_{i}M^{0}_{i}z_{k\hspace{-0.5mm}-\hspace{-0.5mm}d}
\hspace{-1mm}+\hspace{-1mm}\Gamma^{-1}_{i}\sum^{d}_{s=1}(M^{s}_{i}u^{c}_{k\hspace{-0.5mm}-\hspace{-0.5mm}d\hspace{-0.5mm}-\hspace{-0.5mm}s})\Big)\bigg]\Bigg\}.\label{f64}
\end{eqnarray}
In view of the positive definiteness of $\Gamma_{i}$, in order to minimize (\ref{f64}), we take (\ref{f25}) as the optimal controller. Then, the corresponding optimal
cost functional can be expressed as (\ref{f26}). The desired sufficiency is proved.

\bibliographystyle{plain}        
\bibliography{autosam}           

\end{document}